\ifdefined\XeTeXversion\else
\DocumentMetadata{
    lang = en,
    pdfversion = 2.0,
    pdfstandard = ua-2,
    tagging = on,
    tagging-setup = {table/header-columns = 1},
}
\fi

\documentclass[12pt]{article}
\usepackage[a4paper, scale=0.625, hmarginratio=1:1, vmarginratio=1:2, heightrounded]{geometry}

\usepackage{iftex}
\ifPDFTeX
    \usepackage[T1]{fontenc}
    \usepackage{lmodern} 
    \usepackage[utf8]{inputenc}
    \usepackage[english]{babel}
    \usepackage{amssymb}
    \usepackage{mathrsfs}
\else
    \usepackage{polyglossia}
    \usepackage{unicode-math}
    \setdefaultlanguage{english}
    \ifXeTeX
    \else
    \fi
\fi
\usepackage{eucal} 
\usepackage{amsfonts}
\usepackage{amsmath}
\ifPDFTeX\else
    \renewcommand{\boldsymbol}[1]{\symbf{#1}}
\fi
\ifXeTeX
    \ExplSyntaxOn
    \AtBeginDocument
      {
        \NewCommandCopy \kasumVec \Vec
        \RenewDocumentCommand \Vec { m } { \kasumVec { #1 {} } }
        \NewCommandCopy \kasumcheck \check
        \RenewDocumentCommand \check { m }
          { \kasumcheck { \tl_if_head_eq_meaning:nNTF {#1} \mathcal { #1 {} } {#1} } }
      }
    \ExplSyntaxOff
\fi
\usepackage{microtype}
\usepackage[dvipsnames]{xcolor}
\definecolor{introgreen}{RGB}{0,122,0}
\definecolor{urlviolet}{RGB}{128,0,128}

\usepackage{hyperref}
\hypersetup{
    colorlinks=true,
    linktoc=all,
    linkcolor=blue,
    citecolor=introgreen,
    urlcolor=urlviolet,
    pdftitle={Local Definability of HOD in L(R)},
    pdfauthor={Obrad Kasum},
    pdfsubject={Mathematical logic, set theory},
    pdfkeywords={HOD, L(R), inner model theory, determinacy},
}

\newcommand{\doi}[1]{\href{https://doi.org/#1}{\nolinkurl{#1}}}

\usepackage{imakeidx}
\makeindex

\usepackage{amsthm}

\newtheorem*{main}{Main Theorem}
\newtheorem{theorem}{Theorem}[section]
\newtheorem{lemma}[theorem]{Lemma}
\newtheorem{proposition}[theorem]{Proposition}
\newtheorem{corollary}[theorem]{Corollary}

\theoremstyle{definition}
\newtheorem{definition}[theorem]{Definition}
\newtheorem{declaration}[theorem]{Global Assumption}
\newtheorem*{glAss}{Global Assumption}
\newtheorem{notation}[theorem]{Notation}
\newtheorem{remark}[theorem]{Remark}

\newcommand{\case}{\textbf{Case.} }
\newcommand{\claim}{\textbf{Claim.} }

\usepackage{enumitem}

\newlist{parts}{enumerate}{2}
\setlist[parts, 1]{label=\alph*., ref=\alph*}
\setlist[parts, 2]{label=\roman*., ref=\alph{partsi}.\roman{partsii}}

\newlist{assume}{enumerate}{2}
\setlist[assume, 1]{label=\arabic*., ref=\arabic*}
\setlist[assume, 2]{label=\arabic{assumei}.\arabic{assumeii}., ref=\arabic{assumei}.\arabic{assumeii}}

\newlist{pfenum}{enumerate}{4}
\setlist[pfenum, 1]{label=\arabic{pfenumi}$^\circ$}
\setlist[pfenum, 2]{label=\arabic{pfenumii}$'$}
\setlist[pfenum, 3]{label=\arabic{pfenumiii}$''$}
\setlist[pfenum, 4]{label=$\mbox{\arabic*}'''$, leftmargin=*}

\ExplSyntaxOn
\keys_if_exist:nnT { template/list/std } { label }
  {
    \keys_define:nn { template/list/std }
      {
        ref .code:n =
          {
            \tl_set:cn { the \l__block_counter_tl } {#1}
            \__block_normalize_label:c { the \l__block_counter_tl }
          }
      }
  }
\ExplSyntaxOff

\makeatletter
\ExplSyntaxOn
\AtBeginDocument
  {
    \cs_if_exist:NT \__math_tag_dollardollar_display_end:
      {
        \cs_new_eq:NN \__kasum_dollardollar_end: \__math_tag_dollardollar_display_end:
        \cs_gset_protected:Npn \__math_tag_dollardollar_display_end:
          {
            \__kasum_dollardollar_end:
            \int_compare:nNnT { \@listdepth } > { 0 }
              { \tex_parshape:D 1 ~ \@totalleftmargin \linewidth }
          }
      }
  }
\ExplSyntaxOff
\makeatother

\newenvironment{prooff}
{
    \begin{proof}
    \hfill
    \begin{pfenum}
}
{
    \end{pfenum}
    \end{proof}
}

\newcommand{\intro}[1]{\textcolor{introgreen}{#1}\index{#1}}
\newcommand{\dfn}[1]{\textcolor{introgreen}{#1}}
\ExplSyntaxOn
\NewDocumentCommand{\mrg}{m}{\clist_map_inline:nn {#1} {\index{##1}}}
\ExplSyntaxOff

\providecommand{\restriction}{\upharpoonright} 
\newcommand{\chull}{{\mathrm{cHull}}}
\newcommand{\cof}{{\mathrm{cof}}}
\newcommand{\col}{{\mathrm{Col}}}
\newcommand{\cont}{{\mathfrak{c}}}
\newcommand{\crit}{{\mathrm{crit}}}
\newcommand{\dom}{{\mathrm{dom}}}
\newcommand{\hh}{{\mathcal{H}}}
\newcommand{\hod}{{\mathrm{HOD}}}
\newcommand{\hull}{{\mathrm{Hull}}}
\newcommand{\id}{{\mathrm{id}}}
\newcommand{\ii}{{\mathcal{I}}}
\newcommand{\im}{{\mathrm{im}}}
\newcommand{\lh}{{\mathrm{lh}}}
\newcommand{\lp}{{\mathrm{Lp}}}
\newcommand{\M}{{\mathcal{M}}}
\newcommand{\ord}{{\mathrm{Ord}}}
\newcommand{\otp}{{\mathrm{otp}}}
\renewcommand{\P}{{\mathbb{P}}}
\newcommand{\para}{{\parallel}}
\newcommand{\pp}{{\mathcal{P}}}
\newcommand{\proj}{{\downarrow}}
\newcommand{\ps}{{\mathscr{P}}}
\newcommand{\qq}{{\mathcal{Q}}}
\newcommand{\R}{{\mathbb{R}}}
\newcommand{\rank}{{\mathrm{rank}}}
\newcommand{\ran}{{\mathrm{ran}}}
\newcommand{\rest}{{\restriction}}
\newcommand{\T}{{\mathcal{T}}}
\newcommand{\Tau}{{\mathrm{T}}}
\newcommand{\trcl}{{\mathrm{trcl}}}
\newcommand{\ult}{{\mathrm{Ult}}}
\newcommand{\uu}{{\mathcal{U}}}
\newcommand{\ww}{{\mathcal{W}}}

\title{Local Definability of $\hod$ in $L(\R)$}
\author{Obrad Kasum\footnote{
Université Paris Cité, Sorbonne Université, CNRS, IMJ-PRG, F-75013
Paris, France
}}

\begin{document}

\maketitle
\begin{abstract}
    We show that in $L(\R)$, assuming large cardinals in $V$, $\hod\para\eta^{+\hod}$ is locally definable from $\hod\para\eta$ for all $\hod$-cardinals $\eta\in [\boldsymbol{\delta}^2_1,\Theta)$.
    This is a further elaboration of the statement ``$\hod^{L(\R)}$ is a core model below $\Theta$'' made by John Steel.
\end{abstract}
\tableofcontents
\newpage

\section{Introduction}

\subsection{Statement of the Result}

In this paper, we work in a $\mathrm{ZFC}$ universe $V$ with the following large cardinal assumption.

\begin{glAss}
    We assume that there exists a measurable cardinal above a limit of Woodin cardinals.\footnote{
    In a few places, only in this introduction, we are working under opposing assumptions; these are always clearly stated.
    }\qed
\end{glAss}

There are two important consequences of this assumption that will play a central role in the paper.

\begin{theorem}
    The Axiom of Determinacy ($\mathrm{AD}$) holds in $L(\R)$.
\end{theorem}
\begin{proof}
    This is a classical result that came out of the work of Martin, Steel, and Woodin.
    Its full proof can be found in \cite{neeman2010determinacy}.
\end{proof}

The other result concerns the existence and iterability of the mouse $\M_\omega^\sharp$.

\begin{theorem}
    Let $\lambda$ be the least limit of Woodins.
    Then $\M_\omega^\sharp$ exists and it has a unique $(\omega_1,\omega_1+1)$-iteration strategy $\Sigma$.
    The strategy $\Sigma\rest H_{\omega_1}$ is coded by a ${<}\lambda$-universally Baire set.
\end{theorem}
\begin{proof}
    See \ref{47} for a proof sketch.
\end{proof}

We now turn our attention to the model $L(\R)$.
Whenever we write $\Theta$, $\boldsymbol{\delta}^2_1$, $\hod$, and $\hod_x$ (for $x\in\R$), we mean the corresponding objects computed in $L(\R)$.
It was shown by Steel and Woodin that $\hod\para\Theta$\footnote{\label{86}
See Subsec. \ref{Notation} for the use of the symbol $\para$.
} is a premouse.

\begin{theorem}[Steel, Woodin]
    $\hod\para\Theta$ is a domain of a premouse.
\end{theorem}
\begin{proof}
    A proof can be found in \cite{steel2016hod}.
    The authors in fact represent this structure as the direct limit of iterates of an initial segment of $\M_\omega^\sharp$ (with the inherited strategy; see Theorem \ref{1639}).
\end{proof}

The first paper in this direction was \cite{steel1995HOD}, where the author said that ``$\hod^{L(\R)}$ is a core model below $\Theta$''.
However, one usually expects of a core model more than just the structure of a premouse.
Let us briefly recall some of the facts of the core model theory in the $\mathrm{ZFC}$ context.

There is no general definition of what is \textit{the core model}.
Instead, the definition is given on the case-by-case basis, contingent on anti-large-cardinal assumptions.
Of course, the definition of the core model under a weaker anti-large-cardinal assumption is expected to generalize the definition given under a stronger anti-large-cardinal assumption.\footnote{
Cf. \cite[Section 1]{schimmerling2010core}.
}
Let us focus on the anti-large-cardinal assumption ``there is no inner model of $\mathrm{ZFC} + $`there is a Woodin cardinal' ''.
The core model under this assumption was first constructed in \cite{steel1996core}.

Following the source, let us fix a measurable cardinal $\Omega$ and a normal measure $\mu$ on it.
We are interested in studying the universe $V_\Omega$.\footnote{
\cite{steel1996core} develops the core model theory in this scenario.
The assumption that ``$\ord$ is measurable'' was later eliminated -- see \cite{jensen2013k}.
}
The definition of the core model $K$ of $V_\Omega$ is given in \cite[Definition 5.17]{steel1996core}.
Its main properties read as follows.
\begin{theorem}
    Assume $\mathrm{ZFC} + $``$\Omega$ is a measurable cardinal''.
    Suppose there is no inner model of $V_\Omega$ which satisfies $\mathrm{ZFC} + $``there is a Woodin cardinal''.
    Then the following hold.
    \begin{parts}
        \item \label{120'} $K$ is a premouse of height $\Omega$.

        \item \label{122} $K$ is $(\omega, \Omega + 1)$-iterable.

        \item \label{118} For every poset $\P \in V_\Omega$, it holds in $V^\P$ that $K = K^V$.

        \item \label{120} $K \para \omega_1^V$ is $\Sigma_1$-definable over $J_{\omega_1}(\R)$.

        \item \label{128} For $\mu$-almost every $\alpha < \Omega$, $\alpha^{+K} = \alpha^+$.
    \end{parts}
\end{theorem}
\begin{proof}
    The listing is based on \cite[Theorem 1.1]{jensen2013k}.
    The part \ref{120'} is true by definition.
    The parts \ref{122}, \ref{118} and \ref{128} follow from \cite[Theorem 5.18]{steel1996core}.
    The part \ref{120} follows from \cite[Theorem 6.15]{steel1996core}.
\end{proof}

Part \ref{120}, which is called \textit{an inductive definition}, is of particular interest for us.
Namely, the original definition of $K$ is second-order over $V_\Omega$,\footnote{
Cf. \cite[§6, the first paragraph]{steel1996core}.
} so it is not obvious that the initial segment $K \cap H_{\omega_1}$ of $K$ can be determined by only using the limited information of $J_{\omega_1}(\R)$.
In addition, the following is true by putting together parts \ref{118} and \ref{120} of the previous theorem.

\begin{corollary}
    Assume $\mathrm{ZFC} + $``$\Omega$ is a measurable cardinal''.
    Suppose there is no inner model of $\mathrm{ZFC} + $``there is a Woodin cardinal''.
    Let $\kappa$ be an uncountable cardinal.
    Then $K \para \kappa$ is $\Sigma_1$-definable over $L(H_\kappa)$.\footnote{
    Cf. \cite[page 709, paragraph 4]{jensen2013k}.
    }
    \qed
\end{corollary}

This can be even further improved above $\omega_2^V$.

\begin{theorem}\label{144}
    Assume $\mathrm{ZFC} + $``$\Omega$ is a measurable cardinal''.
    Suppose there is no inner model of $\mathrm{ZFC} + $``there is a Woodin cardinal''.
    Let $\alpha \geq \omega_2^V$ be a cardinal of $K$.
    Then $K \para \alpha^{+K}$ is definable (as a predicate) over $H(\alpha^+)$ using only $K \para \alpha$ as a parameter.
\end{theorem}
\begin{proof}
    This follows from \cite[Lemma 3.5]{gitik2006pcf}.
\end{proof}

This establishes existence of a \textit{local} inductive definition of $K$: $K \para \alpha^{+K}$ is determined from $K \para \alpha$ using only the information in $H(\alpha^+)$.
We establish a similar result for the ``core model'' $\hod \para \Theta$ of $L_\Theta(\R)$.
As it was the case with $K$, $\hod$ is defined with reference to the whole background universe $L(\R)$, so it is not a priori obvious that its bounded parts can be recovered using bounded information from the background universe.
Before stating our Main Theorem, let us note that it makes reference to the restriction of the club filter on $[\eta]^\omega$ to $L(\R)$.
We denote this restriction by $\mu_\eta$.
This filter belongs to $L(\R)$ and is a supercompactness measure in that model;\footnote{
See Theorem \ref{1896}.
} in particular, the local structure we exhibit really belongs to $L(\R)$.

\begin{main}
    Suppose $\eta \in [\boldsymbol{\delta}^2_1, \Theta)$ is a cardinal of $\hod$ and $\chi_\eta$ is the least ordinal $\chi > \eta$ such that $L_\chi((\hod \para \eta)^\omega)[\mu_\eta]$ satisfies $\mathrm{ZF}^- + \mbox{``}\ps([\eta]^\omega)\mbox{ exists''}.$\footnote{
    $\chi_\eta$ is well-defined and strictly smaller than $\Theta$ -- see Lemma \ref{2151}.
    } Then $\hod \para \eta^{+\hod}$ belongs to and is definable over $L_{\chi_\eta}((\hod \para \eta)^\omega)[\mu_\eta]$ using the predicates $\in$ and $\mu_\eta$ and the parameter $\hod \para \eta$.
\end{main}
\begin{proof}
    See Corollary \ref{2715}.
\end{proof}

In Theorem \ref{144}, the local structure over which the inductive definition of $K \para \alpha^{+K}$ is given is $H(\alpha^+)$.
The choice of this structure is justified by the stratification
$$V_\Omega = \bigcup_{\omega_2 \leq \kappa < \Omega} H(\kappa^+).$$
Similarly, in our Main Theorem, the local structure is chosen by considering the stratification
$$L_\Theta(\R) = \bigcup_{\boldsymbol{\delta}^2_1 \leq \eta < \chi < \Theta} L_\chi ((\hod \para \eta)^\omega)[\mu_\eta].$$
The underlying point here is that
$$L_\Theta(\R) = L_\Theta ((\hod \para \eta)^\omega)[\mu_\eta],$$
whenever $\eta$ is a $\hod$-cardinal.
This somewhat unusual representation is (at least partially) justified by the description of determined sets beyond $L(\R)$ obtained in \cite{kasum2024building}.





It is not clear how one would get rid of the additional predicate $\mu_\eta$ in the Main Theorem.
This measure seems to play a crucial role in our argument, which will become clear when we give its outline in Subsection \ref{358}.

\subsection{Notation}
\label{Notation}

We review the notation which we will be using.
We will try to follow the notation of \cite{steel2023comparison} as closely as possible and note the differences when they arise.
We will need to use J-hierarchy above some $x$, which is defined as follows:
\begin{eqnarray*}
    J_1^E(x) & := & V_\omega\cup\trcl(\{x\}) \\
    J_{\alpha+1}^E(x) & := & \mathrm{rud}_E(J_\alpha^E(x)\cup\{J_\alpha^E(x)\}) \\
    J_{\gamma}^E(x) & := & \bigcup_{\xi<\gamma}J_\xi^E(x)\quad (\gamma\mbox{ limit}).
\end{eqnarray*}
An $x$-J-structure $M$ has the form
$$M=(J_\alpha^E(x),\in,E,A,\ y : y\in\trcl(\{x\})),$$
where $A$ is amenable to $J_\alpha^E(x)$.
We denote
$$\hat o(M):=\alpha\mbox{ and } o(M):=\omega\alpha.$$
Following \cite{steel1983scales}, we define the first projectum of $M$ as below.

\begin{definition}
\label{projectum}
    Suppose that
    \begin{assume}
        \item $x$ is a set,

        \item $M=(J_\alpha^E(x),\dots)$ is an $x$-J-structure.
    \end{assume}
    Then $\rho_1(M)$ is the least $\rho\leq\alpha$ such that there exists $A\in\boldsymbol\Sigma_1^M$ satisfying 
    $$A\cap J^E_\rho(x)\not\in M.$$
    \qed
\end{definition}

Note that it can (and does) happen, according to our definition, that 
$$\rho_1(M)=1.$$
If $x=0$, the usual definition would have the projectum be $\omega$, so this is a difference to keep in mind.\footnote{
However, note that the notation $\rho_1 = 1$ agrees with the notation of \cite{schimmerling2010core}.
}
If $\rho_1(M)>1$, there are no discrepancies of this kind.
Regarding the notions of the standard parameter and soundness, they have their usual description, but w.r.t.\ the language of $x$-J-structures.
This comment about the language was not necessary when we were defining the projectum: the set $A$ above needs only be boldface definable anyway.
However, in the case of the standard parameter and soundness, it is important to know that the elements of $\trcl(\{x\})$ are always allowed to be used as constants.

All premice are Mitchell-Steel-indexed and they have the soundness degree integrated into them.
If $M$ is a premouse, then $k(M)$ denotes its soundness degree.
For all $(\eta,l)\leq(\hat o(M),k(M))$, we denote
$$M|(\eta,l):=(J_\eta^{E^M},\in,E^M\rest\eta,E^M_\eta,l),$$
$$M|\eta:=M|(\eta,0),$$
$$M\para\eta:=(J_\eta^{E^M},\in, E^M\rest\eta, \emptyset,0).$$
An ordinal $\xi\leq o(M)$ is a \textit{strong cutpoint} of $M$ iff for all extenders $E$ on the $M$-sequence, either $\lh(E)<\xi$ or $\crit(E)>\xi$.
In the case that $M$ has a unique Woodin cardinal, that cardinal is denoted by $\delta(M)$.
All of this notation relativizes to $x$-premice in the obvious way.
If we end up talking about $r\Sigma_{\omega+1}$, the reader should understand this as simply talking about $r\Sigma_\omega$.

If $\T$ is a normal tree of limit length on a premouse $M$, then its \textit{lined-up} part
$$\bigcup_{\xi < \lh(\T)}\M^\T_\xi \para \lh(E^\T_\xi)$$
is denoted by \intro{$\M(\T)$}.
Recall that every premouse occurring on a normal tree extending $\T$ on the position $\lh(\T)$ or later must extend $\M(\T)$.
The height of the lined-up part, $o(\M(\T))$, is denoted by \intro{$\delta(\T)$}.

Let us also point out that we use $\para$ in additional way.
In the $L(\R)$ context, $\hod \para \Theta$ denotes the premouse structure on $\hod \cap V_\Theta$ that comes from Theorem \ref{1639}.
If $\eta < \Theta$, than we interpret $\hod \para \eta$ to mean $(\hod \para \Theta) \para \eta$ (which is then also a premouse).
In this way, whenever the notation $X \para \xi$ is used, the object that it refers to is assumed to have a premouse structure.

\subsection{Prerequisites}
\label{(2)}

The only prerequisite for reading this paper is the basic theory of Mitchell-Steel mice, as presented in \cite{steel2010outline}.
Although not strictly necessary, familiarity with the paper \cite{steel2016hod} is useful: indeed, the starting point for our work is the $\hod$ analysis in $L(\R)$ presented there.
Let us also mention that, even though we are not using results of \cite{steel2023comparison}, our notation is based on the notation established in this book.

The theorem below is a well-known fact in Inner Model Theory, but we were unable to pinpoint an exact reference.
As a courtesy to the reader, we will briefly outline its proof.
Regarding its notation, the mouse $\M_\omega^\sharp$ is defined in \cite[Definition 7.3]{steel2010outline}.
Our convention is that $k(\M_\omega^\sharp)=0$.
This mouse projects to $\omega$ and it is sound.

\begin{theorem}\label{47}
    Let $\lambda$ be the least limit of Woodin cardinals.
    Then $\M_\omega^\sharp$ exists and it has a unique $(\omega_1,\omega_1+1)^*$-iteration strategy $\Sigma$.
    The strategy $\Sigma\rest H_{\omega_1}$ is coded by a ${<}\lambda$-universally Baire set.
    \qed
\end{theorem}
\begin{proof}[Proof outline]
    The existence of $\M_\omega^\sharp$ follows from \cite[Theorem 7.2]{steel2010outline}.
    To verify the uniqueness, we have to consider Q-structures, which will be reviewed in Section \ref{Q-structures}.
    Since $\M_\omega^\sharp$ projects to $\omega$, for all normal trees $\T$ on $\M_\omega^\sharp$ of countable limit length and for all cofinal wellfounded branches $b$ through $\T$, we have that $\qq(\T,b)$ exists.
    If $\T^\frown b$ is according to some $(\omega_1,\omega_1+1)$-iteration strategy, then $\qq(\T,b)$ is $(\omega_1+1)$-iterable.
    This means that any such strategy must pick the same branch through $\T$ (cf.\ \cite[Corollary 6.14]{steel2010outline}).
    This shows the uniqueness of the restriction of $\Sigma$ to normal trees.
    The uniqueness of the full strategy then follows from the normalization (cf.\ \cite{schlutzenberg2021full}).

    To verify that the strategy $\Psi:=\Sigma\rest H_{\omega_1}$ is ${<}\lambda$-universally Baire, fix an uncountable cardinal $\kappa<\lambda$.
    By doing the $K^c$ construction above $\kappa$, we reach $\M_\omega^\sharp$ and obtain the strategy for it from the construction.
    This construction is absolute between $V$ and any generic extension of $V$ by a poset of size ${<}\kappa$.
    In the generic extension, we also reach $\M_\omega^\sharp$ and obtain the strategy for it from the construction.
    Since these strategies come from the realizability into the construction (cf.\ \cite[Theorem 6.6]{steel2010outline}), it is easily seen that club many hulls are ${<}\kappa$-generically correct about $\Psi$.
    This implies that $\Psi$ is ${<}\kappa$-universally Baire (cf.\ \cite[Lemma 4.1]{steel2009derived}).
\end{proof}

The notation $*$ in the statement of the previous theorem means that we are not allowed to make unnecessary drops when stacking normal trees.
By \cite[Theorem 7.2]{steel2010outline}, we can extend the strategy to allow these arbitrary drops as well.
The tradeoff is that we loose uniqueness.
However, we can fix one strategy for $\M_\omega^\sharp$ from the start, so all the arguments will be conditioned on that strategy.

\begin{notation}
    We denote by \intro{$\Sigma_{\M_\omega^\sharp}$} some fixed $(\omega_1, \omega_1 + 1)$-iteration strategy for $\M_\omega^\sharp$.
    \qed
\end{notation}

\subsection{Outline of the Paper}

In Section \ref{Q-structures}, we review the basics of Q-structures.
Under certain circumstances, these structures determine which branches should be picked by an iteration strategy.
They play the crucial role of allowing us to approximate iteration strategies of certain mice inside $L(\R)$.
These mice will be called super-suitable and they will be our primary focus here.
Since $L(\R)$ cannot ascertain the super-suitability, we will need to extract some weaker properties from it.
One such property is the suitability, introduced in Section \ref{Suitable Premice}.
The super-suitability itself is introduced in Section \ref{Super-suitable Premice}, while Section \ref{section: short tree iterability} analyses one a weakening of iterability called \textit{short-tree iterability}.
In Section \ref{HOD as a Direct Limit of Mice}, we describe the already mentioned result of Steel and Woodin on representing $\hod\para\Theta$ as a direct limit of mice.
This concludes the introductory part of the paper.
Sections \ref{2031}-\ref{local definition of hod} contain the main argument of the paper, which is outlined in the following subsection.

\subsection{Outline of the Proof}
\label{358}

Let $\M$ be the premouse $\M_\omega^\sharp$ cut at the successor of its least Woodin cardinal and let $\Sigma_\M$ be its iteration strategy derived from $\Sigma_{\M_\omega^\sharp}$.
Consider the set of all countable, non-dropping, normal iterates of $(\M, \Sigma_\M)$.
A general element of this set is of the form $(N, \Sigma_N)$, where $N$ is a premouse iterate of $\M$ by $\Sigma_\M$ and $\Sigma_N$ is the tail strategy.
This set is naturally equipped with an order: $(N, \Sigma_N)$ is less or equal $(P, \Sigma_P)$ iff $(P, \Sigma_P)$ is a non-dropping normal iterate of $(N, \Sigma_N)$.
In this way, we obtain a directed order.

This directed order indexes a canonical directed family: at the position $(N, \Sigma_N)$, we put the premouse $N$; for $(N, \Sigma_N) \leq (P, \Sigma_P)$, we let the arrow 
$$\pi_{(N, \Sigma_N), (P, \Sigma_P)}$$
be the iteration map from $(N, \Sigma_N)$ to $(P, \Sigma_P)$.
Having obtained a directed family, we can take its direct limit -- call it $\M_\infty$.
It was proved in \cite{steel2016hod} that $\Theta$ is the largest cardinal of $\M_\infty$ and that $\M_\infty \para \Theta = \hod \para \Theta$.

The directed family that we have just described is a surjective image of $\R$.
This means that in $V^{\col(\omega, \R)}$, it becomes countable.
Using its directedness, we can find a sequence cofinal in it.
To pass from an element of this sequence to the next one, we use a non-dropping normal iteration tree.
Therefore, the sequence gives rise to a non-dropping stack of $\omega$ normal trees, and the last model of its branch is $\M_\infty$.
Schlutzenberg proves in \cite{schlutzenberg2021full} that this stack can be normalized, so to obtain a single non-dropping normal tree from $\M$ to $\M_\infty$, and he observes, furthermore, that this tree belongs to $V$.

We will prove that this tree is essentially unique by describing an ``algorithm'' how to to get it.
At this point, we have an explicit description how to get to $\M_\infty$ by a normal tree starting from any premouse in the directed system.
Based on this description, we argue that each of these trees, minus its final branch, actually belongs to $L(\R)$.
To proof proceeds by induction, where the non-trivial step is to show that the branches belong to $L(\R)$.

Our description of the tree identifies its intermediate branches as those which have countably iterable Q-structures.
To prove that one of these branches is in $L(\R)$, we look at the countable hulls of the tree.
In these countable trees, the branches are also picked based on Q-structures.
However, since the hulls are countable, this picking can be done in $L(\R)$.
To get the branch through the big tree, we have to ``put back together'' these countable trees and their branches.
This is done by an ultraproduct argument.

We will explain along the way that for every $\eta < \Theta$, the club filter on $[\eta]^\omega$, when intersected with $L(\R)$, gives rise to a supercompactness measure in $L(\R)$.
This measure is what we use to run the ultrapower argument.
Once we finish proving that the tree minus its final branch is in $L(\R)$, we observe that, in order to get the initial segment of the tree up to some $\eta$, the argument does not really require all $L(\R)$, but rather only the information in $L_{\chi_\eta}((\hod \para \eta)^\omega)[\mu_\eta]$.

To go from here to the local definability, we proceed as follows.
Take any non-dropping normal iterate $M$ of $(\M, \Sigma_\M)$ which has a pre-image for $\eta$ in the direct-limit embedding and let $\T_M$ be the normal tree from $M$ to $\M_\infty$.
Then it is exactly at the stage $\eta$ of the iteration $\T_M$ that we get a model which agrees with $\M_\infty$ up to $\eta$.
However, $\M_\infty$ up to $\eta$ is $\hod$ up to $\eta$, and this is part of the information in the structure $L_{\chi_\eta}((\hod \para \eta)^\omega)[\mu_\eta]$.
In turns out that the $\eta^\mathrm{th}$ model of $\T_M$ actually agrees with $\hod$ up to $\eta^{+\hod}$.
So, here is the local description of $\hod$ up to $\eta^{+\hod}$ based on $\hod$ up to $\eta$: take any non-dropping normal iterate of $(\M, \Sigma_\M)$ that has a pre-image of $\eta$ in the direct-limit embedding, iterate it until you reach $\hod$ up $\eta$, and take the initial segment of that iterate up to its successor of $\eta$ -- this is $\hod$ up to $\eta^{+\hod}$.

The only remaining problem is the fact that $L(\R)$ does not see $\Sigma_\M$ and the direct-limit construction.
This is resolved by working with ``suitable'' premice and short-tree strategies, the same way they did it in \cite{steel2016hod}.

\subsection{AI Disclosure}

AI was \textbf{not} used during the process of writing this paper.
The work had been carried out before AI became useful for mathematical research.

\subsection{Acknowledgments} 

This paper is part of my doctoral dissertation \cite{kasum2025investigation}.
I would like to thank my advisor Grigor Sargsyan for pointing out this problem to me and for many fruitful discussions on the topic.

A part of this paper was written during my stay at the Gdańsk branch of the Institute of Mathematics of the Polish Academy of Sciences.
I am very grateful for their hospitality.

I would also like to thank the anonymous referee for reading the manuscript and making multiple useful remarks.

\subsection{Funding}

I have received funding from the European Union’s Horizon 2020 research and innovation program under the Marie Skłodowska-Curie grant agreement No.\ 945322.

Fondation Sciences Mathématiques de Paris partially supported my visit to the Gdańsk branch of the Institute of Mathematics of the Polish Academy of Sciences, during which part of this paper was written.

\section{Q-structures}
\label{Q-structures}

In this section, we review basic facts about Q-structures.
For an introduction on this, the reader is invited to consult \cite{steel2010outline}.
Our definitions here are based on \cite[Definition 3.3]{mueller2021hod}.

\begin{definition}
    Suppose that
    \begin{assume}
        \item $M$ is a premouse,

        \item $\delta\leq o(M)$,

        \item $Q\unlhd M$.
    \end{assume}
    Then $Q$ is a \intro{Q-structure of $M$ at $\delta$} iff both \ref{a} and \ref{b} hold, where:
    \begin{parts}
        \item\label{a} either $\delta=o(Q)$ or $Q\models$``$\delta$ is Woodin'';

        \item\label{b} one of the conditions \ref{bi}, \ref{bii}, or \ref{biii} is met, where:
        \begin{parts}[ref=\roman{partsii}]
            \item\label{bi} $\hat o(Q)<\hat o(M)$, $k(Q)=\omega$, and $\delta$ is not Woodin in $J_1(M)$\footnote{Note that the universe of $J_1(M)$ is the same as the universe of $M\para(\hat o(M)+1)$.};

            \item\label{bii} $\hat o(Q)=\hat o(M)$, $\rho(Q)=\delta$, and there exists an $r\Sigma_{k(M)+1}$ subset of $\delta$ witnessing that $\delta$ is not Woodin in $J_1(M)$;

            \item\label{biii} $\hat o(Q)=\hat o(M)$ and $\rho(Q)<\delta$.\qed
        \end{parts}
    \end{parts}
\end{definition}

\begin{notation}
    In the setup of the previous definition, there there exists at most one $Q$ which is a Q-structure for $M$ at $\delta$.
    If such $Q$ does exist, we call it \textit{the} Q-structure for $M$ at $\delta$ and we denote it by \intro{$\qq(M,\delta)$}.
    To say that such $Q$ exists, we use the shorthand of saying ``\dfn{$\qq(M,\delta)$ exists}''.\qed 
\end{notation}







In an iteration tree of limit length which is not too complicated, cofinal branches come with a naturally assigned Q-structures.
These structures can then be used to pick the right branch through that tree.

\begin{definition}
    Suppose that
    \begin{assume}
        \item $M$ is a premouse,

        \item $\T$ is a normal tree on $M$ of limit length,

        \item $b$ is a cofinal well-founded branch of $\T$.
    \end{assume}
    Then \intro{$\qq(\T,b)$} exists iff $\qq(\M^\T_b,\delta(\T))$ exists.
    In that case, we denote by $\qq(\T,b)$ the structure $\qq(\M^\T_b,\delta(\T))$.\qed
\end{definition}

The Q-structures corresponding to two different cofinal branches are mutually in comparable, unless they are of a certain particular type.
We isolate this exception in the following definition.

\begin{definition}
    Suppose that
    \begin{assume}
        \item $M$ is a premouse,

        \item $\T$ is a normal tree on $M$ of limit length,

        \item $b$ is a cofinal well-founded branch of $\T$.
    \end{assume}
    Then $(M,\T,b)$ is an \intro{anomaly} iff all of the following conditions are simultaneously met:
    \begin{parts}
        \item $M$ is not sound,

        \item $b$ does not drop,

        \item $\qq(\T,b)$ exists and is equal to $\M^\T_b$.\qed 
    \end{parts}
\end{definition}

\begin{proposition}
\label{252}
    Suppose that
    \begin{assume}
        \item $M$ is a premouse,
        
        \item $\T$ is a normal tree of limit length on $M$,

        \item $b\not=c$ are cofinal well-founded branches through $\T$,

        \item $(M,\T,b)$ and $(M,\T,c)$ are not anomalies.
    \end{assume}
    Then neither $\qq(\T,b)\unlhd\qq(\T,c)$ nor $\qq(\T,c)\unlhd\qq(\T,b)$.
\end{proposition}
\begin{proof}
    See \cite[Theorem 6.12]{steel2010outline}.
\end{proof}

If $\T$ is a normal tree on $M$ of limit length and if the next branch to be picked has the Q-structure, we might try to guess that structure before actually knowing the branch.
We now work towards introducing this structure and finally succeed in doing so in Definition \ref{288}.

\begin{definition}
    Suppose that $M$ is a premouse.
    Then $M$ is \intro{countably iterable} iff for all countable premice $\bar M$ and all elementary $i:\bar M\to M$, it holds that $\bar M$ is $(\omega_1+1)$-iterable.\qed
\end{definition}

\begin{lemma}\label{unique Q-structure}
    Suppose that $P$ is a premouse.
    Then there exists at most one premouse $Q$ satisfying that
    \begin{parts}
        \item $P\unlhd Q$,

        \item $o(P)$ is a strong cutpoint in $Q$,

        \item $Q$ is the Q-structure of $Q$ at $o(P)$,\footnote{Or in other words, $Q=\qq(Q,o(P))$.}

        \item $Q$ is sound above $o(P)$,

        \item $Q$ is countably iterable.
    \end{parts}
\end{lemma}
\begin{proof}
    Let us assume otherwise and let $Q_0\not=Q_1$ be two witnesses.
    Since we can always take a countable hull of some $H_\theta$, for a large enough $\theta$, we may assume w.l.o.g. that $P$, $Q_0$, and $Q_1$ are countable.
    By the argument of \cite[Corollary 3.12]{steel2010outline}, we have that either $Q_0\unlhd Q_1$ or $Q_1\unlhd Q_0$.
    However, the minimality which is a part of the definition of a Q-structure would then imply $Q_0=Q_1$, which is a contradiction.
\end{proof}

\begin{definition}\label{Q-structure above}
    Suppose that $P$ is a premouse.
    \dfn{The Q-structure above $P$}, denoted by \intro{$\qq(P)$}, is the unique premouse $Q$ satisfying:
    \begin{parts}
        \item $P\unlhd Q$,

        \item $o(P)$ is a strong cutpoint in $Q$,

        \item $Q$ is the Q-structure of $Q$ at $o(P)$,

        \item $Q$ is sound above $o(P)$,

        \item $Q$ is countably iterable.\qed
    \end{parts}
\end{definition}

\begin{definition}\label{288}
    Suppose that
    \begin{assume}
        \item $M$ is a premouse,

        \item $\T$ is a normal tree on $M$ of limit length.
    \end{assume}
    Then \dfn{$\qq(\T)$ exists} iff $\qq(\M(\T))$ exists.\mrg{$\qq(\T)$}
    In that case, we define 
    $$\qq(\T):=\qq(\M(\T)).$$
    \qed
\end{definition}

If mice do not have extenders overlapping local Woodins, they are said to be tame.
Tame mice are simple enough so that Q-structures can be used to identify the right branches through trees on them, i.e. the branches that must be pick by any sufficiently strong strategy for those mice.

\begin{definition}
    Suppose that $M$ is a premouse.
    Then $M$ is \intro{tame} iff for all $\eta<\hat o(M)$, if $E^M_\eta\not=\emptyset$, then for all $\delta\in [\crit(E),\eta)$, $M\para\eta\models$``$\delta$ is not Woodin''.\qed
\end{definition}

\begin{proposition}
    Suppose that
    \begin{assume}
        \item $M$ is a tame premouse,

        \item $\T$ is a normal tree on $M$ of limit length.
    \end{assume}
    Then there exists at most one cofinal wellfounded branch $b$ through $\T$ such that
    \begin{parts}
        \item $(M,\T,b)$ is not an anomaly,

        \item $\qq(\T,b)$ exists,

        \item $\qq(\T,b)$ is countably iterable.
    \end{parts}
\end{proposition}
\begin{proof}
    Let us assume otherwise and let $b_0\neq b_1$ be two such branches.
    Since $M$ is tame, $\delta(\T)$ is a strong cutpoint in both $\qq(\T,b_0)$ and $\qq(\T,b_1)$.
    By the uniqueness of Q-structures (cf. Lemma \ref{unique Q-structure}), we have that $\qq(\T,b_0)=\qq(\T,b_1)$, which contradicts Proposition \ref{252}.
\end{proof}







\begin{proposition}
    Suppose that
    \begin{assume}
        \item $M$ is a tame premouse,

        \item $\Sigma_M$ is an $(\omega_1+1)$-iteration strategy for $M$,
        
        \item $\T$ is a normal tree on $M$ of limit length according to $\Sigma_M$,

        
        \item $b:=\Sigma_M(\T)$,

        \item $\qq(\T,b)$ exists.
    \end{assume}
    Then $\qq(\T)$ exists and is equal to $\qq(\T,b)$.
\end{proposition}
\begin{proof}
    The Q-structure $\qq(\T,b)$ is countably iterable since $\M^\T_b$ is countably iterable.
    The ordinal $\delta(\T)$ is a strong cutpoint in $\qq(T,b)$ because $M$ is tame.
    This suffices for the conclusion.
\end{proof}

The ordinal $(\boldsymbol{\delta}^2_1)^{L(\R)}$ will figure prominently in the present work.
There are many different characterizations of this ordinal, but we choose the one most useful for our purposes (cf. \cite[Lemma 1.12]{steel1983scales}).
Since we will not compute this ordinal in any other model except $L(\R)$, we omit the superscript.

\begin{definition}\label{Sigma-1 reflection}
    The ordinal $\boldsymbol{\delta}^2_1$ is the least ordinal $\delta$ satisfying that $\Sigma_1$-formulas with parameters in $\R\cup\{\R\}$ are absolute between $L_\delta(\R)$ and $L(\R)$.\qed
\end{definition}

This reflection implies that all countable mice that have $\omega_1$-iteration strategies in $L(\R)$, have such strategies in $L_{\boldsymbol{\delta}^2_1}(\R)$.
The following proposition is an example of how this fact can be used.

\begin{proposition}\label{412}\label{Q-structures L(R)}
    Suppose that
    \begin{assume}
        
        \item $M$ is a countable $\omega$-small\footnote{See \cite[Definition 7.1]{steel2010outline}} premouse,

        \item for all $\delta<o(M)$, if $M\models$``$\delta$ is Woodin'', then $\qq(M,\delta)$ exists,

        \item $M$ is $(\omega_1+1)$-iterable.
    \end{assume}
    Then $L(\R)\models$``$M$ is $(\omega_1+1)$-iterable''.
    Moreover, there exists a set of reals in $L_{\boldsymbol{\delta}^2_1}(\R)$ which canonically codes an $\omega_1$-iteration strategy for $M$.
\end{proposition}
\begin{proof}
    Since $M$ is $(\omega_1+1)$-iterable, we have that it is weakly $\omega$-iterable\footnote{See \cite[Definition 7.7]{steel2010outline}}.
    By the proof of \cite[Theorem 7.10]{steel2010outline}, we have that $M$ has an $(\omega_1+1)$-iteration strategy in $L(\R)$.
    The moreover part follows from Definition \ref{Sigma-1 reflection}.
\end{proof}


\section{Suitable Premice}
\label{Suitable Premice}

We are really interested in super-suitable premice.
They are the appropriate initial segments of iterates of $\M_\omega$ and they are going to be introduced in Section \ref{Super-suitable Premice}.
However, super-suitable premice cannot be defined internally in $L(\R)$, so we need to work with an approximate notion, that of a suitable premouse.
Of course, it will be a theorem that super-suitable premice are suitable.

\begin{definition}
    Suppose that
    \begin{assume}
        \item $a$ is countable and transitive,

        \item $M$ is an $a$-premouse.
    \end{assume}
    Then $M$ is \intro{$\lp$-good} iff all of the following conditions are met:
    \begin{parts}
        \item $k(M)=\omega$,
        
        \item $\rho_\omega(M)=1$\footnote{See Definition \ref{projectum} and the comment after it.},
        
        \item $M$ has an $\omega_1$-iteration strategy in $L(\R)$.\qed
    \end{parts}
\end{definition}

$\lp$-good $a$-premice extend each other and there is no longest one among them.
We will be interested in their supremum.

\begin{lemma}
    Suppose that
    \begin{assume}
        \item $a$ is countable and transitive,



        \item $M,N$ are $\lp$-good.
    \end{assume}
    Then either $M\unlhd N$ or $N\unlhd M$.
\end{lemma}
\begin{proof}
    This is a straightforward generalization of \cite[Corollary 3.12]{steel2010outline}.
\end{proof}

\begin{lemma}
    Suppose that $a$ is countable and transitive.
    Then for all $\lp$-good $a$-premice $M$, there exists an $\lp$-good $a$-premouse $N$ such that $M\lhd N$.
\end{lemma}
\begin{proof}
    Since $M$ projects to $1$ and is sound, there exists a surjection $f:J_1(a)\twoheadrightarrow M$ which is definable over $M$.
    Now, if look at $J_1(M)$, organized as an $a$-premouse, we see that it is obtained as the rudimentary closer of $M\cup\{M\}$.
    However, the rudimentary functions can be listed recursively, so the facts that $V_\omega\cup\{f\}\subseteq J_1(a)$ allow us to define a surjection
    $$g:J_1(a)\twoheadrightarrow J_1(M)$$
    over $J_1(M)$.
    This means that $J_1(M)$ projects to $1$ and is consequently $\lp$-good.
\end{proof}

Supremum of all $\lp$-good $a$-mice is denoted by $\lp(a)$.
This object is itself an $a$-mouse and it looks like the power set of $a$.
Since sets appearing in $\lp(a)$ are not too complicated, i.e. they come from mice that have strategies in $L(\R)$, this object can be understood as a \textit{lower part} of the full powerset of $a$ (hence the abbreviation $\lp$).

\begin{definition}\label{lp(a)}
    Suppose that $a$ is countable and transitive.
    Then \intro{$\lp(a)$} is the unique $a$-premouse satisfying:
    \begin{parts}
        \item\label{390} for all $\lp$-good $a$-premice $M$, we have that $M\unlhd\lp (a)$;

        \item no proper initial segment of $\lp(a)$ satisfies the previous close.\qed
    \end{parts}
\end{definition}

\begin{lemma}
\label{487}
    Suppose that $a$ is countable and transitive.
    Then the following holds:
    \begin{parts}
        \item\label{401} $k(\lp(a))=0$,

        \item\label{403} $\lp(a)$ is countable.
    \end{parts}
\end{lemma}
\begin{prooff}
    \item If it were the case that $k(\lp(a))>0$, then the $a$-premouse obtained from $\lp(a)$ by decreasing $k(\lp(a))$ by 1 would also satisfy condition \ref{390} of Definition \ref{lp(a)}, while being a strict initial segment of $\lp(a)$.
    This shows that \ref{401} must hold.
    
    \item To establish \ref{403}, let us assume otherwise.
    Then the set $C$ of all $\alpha\in (\rank(a),\omega_1)$ such that there exists an $\lp$-good $a$-premouse $M$ satisfying $\hat o(M)=\alpha$ is cofinal in $\omega_1$.

    \item For all $\alpha\in C$, an $M$ witnessing this fact is unique and we denote it by $M_\alpha$.
    
    \item Let $(f_\alpha : \alpha\in C),g$ be as follows:
    \begin{assume}
        \item for all $\alpha\in C$, $f_\alpha:a\cup\{a\}\twoheadrightarrow M_\alpha$ is given by the soundness,

        \item $g:\omega\twoheadrightarrow a\cup\{a\}$ is an arbitrary enumeration,

        \item for all $\alpha\in C$, $e_\alpha:=\{(m,n) : f_\alpha(g(m))\in f_\alpha(g(n))\}\subseteq\omega^2$.
    \end{assume}
    
    \item The sequence $(e_\alpha : \alpha\in C)$ is injective and it belongs to $L(\R)$.
    This contradicts $L(\R)\models\mathrm{AD}$.
\end{prooff}

We said that $\lp(a)$ should be understood as a lower part of the powerset of $a$.
It turns out that if $a$ is countable, transitive, and self-wellorderable\footnote{A transitive set $a$ is \textit{self-wellorderable} iff $J_1(a)\models$``there exists a wellordering on $a$''.}, it is in fact the case that
$$\lp(a)=H(|a|^+)^{\M_\omega (a)}=H(|a|^+)^{\hod^{L(\R)}(a\cup\{a\})}.$$
This is the content of \cite[Theorem 6.4]{steel2016hod}, which we reproduce here.

\begin{theorem}
\label{484}\label{(6)t}
    Suppose that
    \begin{assume}
        \item $a$ is countable and transitive,

        \item $b\subseteq a$.
    \end{assume}
    Then the following are equivalent.
    \begin{parts}
        \item $b\in\lp(a)$.

        \item $b\in\M_\omega(a)$.
        
        \item $b$ is definable over $(L(\R),\in)$ from parameters in $\ord\cup a\cup\{a\}$.\qed 
    \end{parts}
\end{theorem}

\begin{corollary}
    Suppose that $a$ is countable, transitive, and self-wellorderable.
    Then it holds that
    \begin{parts}
        \item\label{458} $\lp(a)=H(|a|^+)^{\M_\omega(a)} = H(|a|^+)^{\hod^{L(\R)}(a\cup\{a\})}$,
        
        \item\label{460} $\lp(a)\models\mathrm{ZFC}^-$,

        \item\label{462} $\lp(a)$ is $\omega$-sound.\qed

    \end{parts}
\end{corollary}

We are ready to introduce the notion of a suitable premouse.
Here and later, we will need reorganize premice into premice over their initial segments (when this is possible). 

\begin{notation}
    Suppose that $M$ is a premouse and $\delta$ is a strong cutpoint of $M$.
    Then we denote by \intro{$M/\delta$} the canonical reorganization of $M$ into an $(M\para\delta)$-premouse.\qed 
\end{notation}

\begin{definition}
    Suppose that $M$ is a premouse.
    Then $M$ is \intro{suitable} iff there exists $\delta<\omega_1$ such that
    \begin{parts}
        
        \item $M\models\mathrm{ZFC}^-+$``$\delta$ is the largest cardinal''$+$``$\delta$ is Woodin'',

        \item $\delta$ is a strong cutpoint of $M$,

        \item for all $\eta\leq\delta$, $\lp(M\para\eta)\subseteq M$,

        \item for all $\eta\leq\delta$, if $\eta$ is a strong cutpoint of $M$, then $\lp (M\para\eta)=(M\para\eta^{+M})/\eta$,

        \item for all $\eta<\delta$, $\lp(M\para\eta)\models$``$\eta$ is not Woodin''.\qed
    \end{parts}
\end{definition}

Note that in order for $M$ to be suitable, one has to be able to find a \textit{countable} ordinal $\delta$ witnessing this fact.
An immediate consequence is the fact that all suitable premice are countable.

\begin{lemma}
    Suppose that $M$ is suitable.
    Then $M$ is countable.
\end{lemma}
\begin{proof}
    This follows from Lemma \ref{487}. and the fact that there exists $\delta<\omega_1$ such that $M/\delta=\lp(M\para\delta)$.
\end{proof}


\section{Super-suitable Premice}
\label{Super-suitable Premice}

A super-suitable premouse is obtained from an iterate of $\M_\omega$ by cutting it at the successor of its least Woodin.
We make precise in next few definitions.
Here and later, we shall introduce several variations on the notion of an iterate.
What we call here simply ``an iterate'' is elsewhere called (more cumbersomely) ``a nondropping iterate''.
Since we will not have the need to talk about dropping iterates, we omit this additional qualifier.
Similarly, since we will not have the need to talk about uncountable iterates, we omit the qualifier ``countable'' as well and incorporate the countability in the definitions.

\begin{definition} 
    Suppose that 
    \begin{assume}
        \item  $M,N$ are countable premice,

        \item $\Sigma_M$ is an $(\omega_1+1)$-iteration strategy for $M$.
    \end{assume}
    Then $N$ is a \intro{normal $\Sigma_M$-iterate} of $M$ iff there exists a countable normal tree $\T$ on $M$ according to $\Sigma_M$ whose last model is $N$ and whose main branch does not drop.\qed
\end{definition}

\begin{definition} 
    Suppose that 
    \begin{assume}
        \item  $M,N$ are countable premice,

        \item $\Sigma_M$ is an $(\omega_1,\omega_1+1)$-iteration strategy for $M$.
    \end{assume}
    Then $N$ is a \intro{$\Sigma_M$-iterate} of $M$ iff there exists a countable stack of countable normal trees on $M$ according to $\Sigma_M$ whose last model is $N$ and whose main branch does not drop.\qed
\end{definition}

\begin{definition} \label{724}
    Suppose that 
    \begin{assume}
        \item $\tau:=\delta(\M_\omega)^{+\M_\omega}$,

        \item $\M:=\M_\omega\para\tau$,

        \item $\Sigma_\M$ is the iteration strategy for $\M$ obtained from $\Sigma_{\M_\omega^\sharp}$,

        \item $M$ are countable premice.
    \end{assume}
    Then we define the following.
    \begin{parts}
        \item $M$ is \intro{super-suitable} iff $M$ is a $\Sigma_\M$-iterate of $\M$.

        \item If $M$ is super-suitable, then \intro{$\Sigma_M$} is the $(\omega_1,\omega_1+1)$-iteration strategy for $M$ induced by $\Sigma_\M$.
        \qed
    \end{parts}
\end{definition}

Since the supremum of Woodin cardinals of $\M_\omega$ is countable and since all iterates are countable by our choice of the definition, we have that all super-suitable mice are countable.
Furthermore, since all super-suitable mice are elementarily equivalent to the mouse $\M$ of the previous definition, they are all tame.
We highlight this in the following lemma.

\begin{lemma}
    All super-suitable premice are countable and tame.\qed
\end{lemma}

When we consider super-suitable mice, we will only consider them together with their canonical strategies.
Thus, we can simplify the terminology of iterates a bit.

\begin{notation}
    Suppose that $M$ is super-suitable and that $N$ is a countable premouse.
    Then we define the following.
    \begin{parts}
        \item $N$ is a \intro{normal iterate} of $M$ iff $N$ is a normal $\Sigma_M$-iterate of $M$.

        \item $N$ is an \intro{iterate} of $M$ iff $N$ is a $\Sigma_M$-iterate of $M$.\qed
    \end{parts}
\end{notation}

As we have already pointed it out, the suitability is an approximation to the super-suitability, so the following proposition is to be expected.

\begin{proposition}\label{536}
    Suppose that $M$ is a super-suitable premouse.
    Then $M$ is suitable.
\end{proposition}
\begin{proof}
    See the paragraph immediately below Definition 6.8 of \cite{steel2016hod}.
\end{proof}

If $M$ is super-suitable, then $M$ is in particular tame.
This has for a consequence that $\delta(\T)$ is a strong cutpoint in $\qq(\T,b)$ whenever $\T$ is a normal tree on $M$, $b$ a cofinal wellfounded branch through $\T$, and $\qq(\T,b)$ exists.
In particular, the $\M(\T)$-premouse $\qq(\T,b)/\delta(\T)$ is defined.
If $\qq(\T,b)$ is iterable, then the strategy $\Sigma_M$ must pick the branch $b$ for the tree $\T$, which is the content of the following proposition.

\begin{proposition}\label{canonical strategies, Q-structures} 
    Suppose that
    \begin{assume}
        \item $M$ is a super-suitable premouse,

        \item $\T$ is a countable normal tree on $M$ of limit length according to $\Sigma_M$,

        \item $b$ is a cofinal wellfounded branch through $\T$,

        \item $\qq(\T,b)$ exists and is $(\omega_1+1)$-iterable.
    \end{assume}
    Then $\qq(\T,b)/\delta(\T)\unlhd\lp(\M(\T))$ and $\Sigma_M(\T)=b$.
\end{proposition}
\begin{prooff}
    \item\label{570} By Proposition \ref{412}, we have that $\qq(\T,b)$ has an $(\omega_1+1)$-iteration strategy in $L(\R)$.
    In particular, $\qq(\T,b)$ is $\lp$-good and consequently, 
    $$\qq(\T,b)/\delta(\T)\unlhd\lp(\M(\T)).$$
    
    \item Let us now verify that $\Sigma_M(\T)=b$.
    We assume towards contradiction that the branch $c:=\Sigma_M(\T)$ is distinct from $b$.

    \item\label{664} The Q-structure $\qq(\T,c)$ does not exist, for otherwise it would be $(\omega_1+1)$-iterable, while there can be only one branch with an iterable Q-structure.
    Hence, we get that $c$ does not drop and $\M^\T_c\models$``$\delta(\T)$ is Woodin''.

    \item\label{667} Since $c$ does not drop, we have that $\M^\T_c$ is super-suitable.
    
    \item Putting together \ref{664} and \ref{667}, we get that $\M^\T_c$ is suitable and that $\delta(\T)$ is its unique Woodin.
    This means that
    $$\M^\T_c/\delta(\T)=\lp (\M^\T_c\para\delta(\T))=\lp (\M(\T)).$$

    \item Adding to this the conclusion of \ref{570}, we have that 
    $$\qq(\T,b)\unlhd\lp (\M(\T))=\M_c^\T.$$
    
    \item Since $\delta(\T)$ is Woodin in $\M_c^\T$, it follows that $\qq(\T,b)=\M_c^\T$.

    \item Now, this would mean that $\qq(\T,c)$ does exist (and is equal to $\qq(\T,b)$), which is in contradiction with \ref{664}.
\end{prooff}










\section{Short tree Iterability}
\label{section: short tree iterability}

The model $L(\R)$ does not contain an $\omega_1$-itertion strategy for $\M_\omega$.
The reason for this is that otherwise one could easily construct, using the genericity iterations, $\R^\sharp$ inside $L(\R)$.
A consequence of this fact is that $L(\R)$ cannot check if a given countable normal tree on the suitable initial segment of $\M_\omega$ is according to the strategy for $\M_\omega$.
This also means that one cannot expect to be able to verify inside $L(\R)$ whether a premouse is super-suitable and for this reason, we introduced the notion of suitability.
Now, given a super-suitable premouse $M$, we can verify inside $L(\R)$ that it is suitable, but we cannot know the strategy $\Sigma_M$.
What we do have is a partial strategy, one which knows how to continue the tree as long as the tree is \textit{short} and breaks down once the tree stops being short (such trees are called \textit{maximal}).
This partial strategy is called \textit{the short tree strategy} and denoted by $\Sigma^\mathrm{st}_M$.


\begin{definition}
    Suppose that
    \begin{assume}
        \item $M$ is suitable,

        \item $\T$ is a countable normal tree on $M$ of limit length.
    \end{assume}
    The predicates ``\dfn{$\T$ is short}'', ``\dfn{$\T$ is maximal}'', ``\dfn{$\T\in\dom(\Sigma^\mathrm{st}_M)$}'' and the value \dfn{$\Sigma^\mathrm{st}_M(\T)$} are defined by recursion on $\lh(\T)$, as follows.\mrg{short, maximal, $\Sigma^\mathrm{st}_M$}
    If $\T$ is according to $\Sigma^\mathrm{st}_M$, then
    \begin{parts}
        \item $\T$ is short iff $\qq(\T)$ exists,

        \item $\T$ is maximal iff it is not short,

        \item $\T\in\dom(\Sigma^\mathrm{st}_M)$ iff $\T$ is short and there exists a cofinal wellfounded branch $b$ through $\T$ such that $\qq(\T)\unlhd\M^\T_b$,

        \item if $\T\in\dom(\Sigma^\mathrm{st}_M)$, the $\Sigma^\mathrm{st}_M(\T)$ is the unique cofinal wellfounded branch $b$ through $\T$ satisfying $\qq(\T)\unlhd\M^\T_b$.\qed
    \end{parts}
\end{definition}

To explain the previous definition, at successor stages, the tree is built by picking and applying an extender, so the strategy is concerned by limit stages only.
If the tree is built according to the short tree strategy and if we reach a limit stage, several things can occur.
One possibility is that the tree is maximal, in which case the strategy breaks down.
In particular, no extension of such a tree will be according to the short tree strategy.
The other possibility is that the tree is short.
In this case, we would like to say that $\Sigma^\mathrm{st}_M(\T)$ is defined, but the fact that $\qq(\T)$ exists does not guarantee that there is a branch whose Q-structure is exactly $\qq(\T)$.
This last scenario is a pathology and we would like to exclude it.
Assuming that this pathology does not occur, we then know that the branch $b$ whose Q-structure is equal to $\qq(\T)$ is unique, so we can set $\Sigma^\mathrm{st}_M(\T):=b$.
The pathology just mentioned is dealt with by introducing the notion of short tree iterability.
This notion also excludes another pathology, that of the possibility that an ultrapower in the tree may produce an illfounded model, and it additionally guarantees that the iterates of $M$ stay suitable (which is sometimes called \textit{being fullness preserving}).

\begin{definition}\label{short tree iterable}
    Suppose that $M$ is suitable.
    Then $M$ is \intro{short tree iterable} iff for all countable trees $\T$ on $M$ according to $\Sigma^\mathrm{st}_M$, we have that
    \begin{parts}
        \item if $\T$ has a last model $N$, then 
        \begin{parts}
            \item $\T$ can be normally extended by any extender on the $N$-sequence without producing an illfounded model,

            \item if the branch $M$-to-$N$ of $\T$ does not drop, then $N$ is suitable,
        \end{parts}

        \item if $\T$ is short, then $\T\in\dom(\Sigma^\mathrm{st}_M)$,

        \item if $\T$ is maximal, then there exists a nondropping cofinal wellfounded branch $b$ through $\T$ such that $\M_b^\T$ is suitable.\qed
    \end{parts}
\end{definition}














We have already seen that a super-suitable premouse is suitable.
In this sense, the suitability is an approximation to the super-suitability.
We want to see next that $\Sigma^\mathrm{st}_M$ is an approximation to $\Sigma_M$ whenever $M$ is super-suitable.

\begin{proposition}\label{strategy and st strategy}
    Suppose that $M$ is super-suitable.
    Then $\Sigma_M^\mathrm{st}\subseteq\Sigma_M$.
\end{proposition}
\begin{prooff}
    \item By induction on $\lh(\T)$, we show that for all $\T\in\dom(\Sigma_M^\mathrm{st})$, if
    $$\forall\eta<\lh(\T)\mbox{ limit},\ \T\rest\eta\in\dom(\Sigma_M),$$
    then $\T\in\dom(\Sigma_M)$ and $\Sigma_M(\T)=\Sigma_M^\mathrm{st}(\T)$.
    
    \item Since $\T$ is according to $\Sigma_M$, we immediately have that $\T\in\dom(\Sigma_M)$.

    \item Let us denote by $b$ the branch $\Sigma_M^\mathrm{st}(\T)$.
    By the definition of $\Sigma_M^\mathrm{st}$, we have that $\qq(\T)$ exists and satisfies $\qq(\T)\unlhd\M^\T_b$.
    Thus, $\qq(\T,b)=\qq(\T)$ exists and is $(\omega_1+1)$-iterable.

    \item By Proposition \ref{canonical strategies, Q-structures}, we must have that $\Sigma_M(\T)=b$, as required.
\end{prooff}

Super-suitable mice are not only suitable, but also iterable.
Since we are interested in $L(\R)$-approximations, we should verify the short tree iterability of these mice.

\begin{proposition}\label{super-suitable, short tree iterable}
    Suppose that $M$ is super-suitable.
    Then $M$ is short tree iterable.
\end{proposition}
\begin{prooff}
    \item Let $\T$ be a countable normal tree on $M$ according to $\Sigma^\mathrm{st}_M$.
    By Proposition \ref{strategy and st strategy}, $\T$ is also according to $\Sigma_M$.

    \item What needs to be verified depends on the type of the tree $\T$, as can be seen from Definition \ref{short tree iterable}.
    We have three cases.

    \item\case $\T$ is of successor length.
    \begin{proof}
        First, we need to verify that when $\T$ is extended by one extender, we do not get an illfounded model.
        However, this is obvious since $\T$ is according to $\Sigma_M$ and $\Sigma_M$ is an $(\omega_1,\omega_1+1)$-iteration strategy for $M$.
        The second thing to verify is that if the main branch of $\T$ does not drop, then the last model of $\T$ is suitable.
        For this, observe that if the main branch of $\T$ does not drop, then the last model of $\T$ is an iterate of a super-suitable premouse $M$, so it must itself be super-suitable.
    \end{proof}

    \item\case $\T$ is short.
    \begin{prooff}
        \item Being short means that $\qq(\T)$ exists, so we only need to verify that there exists a cofinal wellfounded branch $b$ through $\T$ such that $\qq(\T)\unlhd\M^\T_b$ works.
        We claim that $b:=\Sigma_M(\T)$ works.

        \item Let us assume towards contradiction that $b$ does not drop.
        Then $\M^\T_b$ is suitable and $\delta(\T)$ is its Woodin.

        \item By definition, the premouse $\qq(\T)$ is $(\omega_1+1)$-iterable and $\delta(\T)$ is its strong cutpoint.
        This means that $\qq(\T)/\delta(\T)$ is $\lp$-good and consequently,
        $$\qq(\T)/\delta(\T)\unlhd\lp(\M(\T)).$$

        \item However, we have by suitability that 
        $$\M^\T_b/\delta(\T)=\lp(\M^\T_b\para\delta(\T))=\lp(\M(\T)),$$
        so we conclude that $\qq(\T)\unlhd\M^\T_b$.

        \item If it were the case that $o(\qq(\T))< o(\M^\T_b)$, it would hold that $\delta(\T)$ is not Woodin in $\M^\T_b$.
        This being absurd, we conclude that the J-structures underlying premice $\qq(\T)$ and $\M^\T_b$ are equal.

        \item However, this is a contradiction since $\rho(\qq(\T))\leq\delta(\T)$ while 
        $$\rho_\omega(\M^\T_b)=\hat o(\M^\T_b).$$

        \item The conclusion is then that $b$ drops, which means that $\qq(\T,b)$ exists.

        \item By the uniqueness of the Q-structure above $\M(\T)$ (cf. Proposition \ref{unique Q-structure}), it follows that $\qq(\T,b)=\qq(\T)$ and consequently, $\qq(\T)\unlhd\M^\T_b$, as required.
    \end{prooff}

    \item\case $\T$ is maximal.
    \begin{prooff}
        \item In this case, we need to verify that there exists a nondropping cofinal wellfounded branch $b$ through $\T$ such that $\M^\T_b$ is suitable.
        We claim that $b:=\Sigma_M(\T)$ works.

        \item If $b$ dropped, then $\qq(\T,b)$ would exist and be $(\omega_1+1)$-iterable.
        In other words, $\qq(\T)$ would exist and be equal to $\qq(\T,b)$, which is not the case.

        \item Thus, $b$ does not drop.
        
        \item It follows that $\M^\T_b$ is an iterate of a super-suitable premouse $M$, as witnessed by $\T^\frown b$, so $\M^\T_b$ is super-suitable.

        \item In particular, $\M^\T_b$ is suitable, as required.
    \end{prooff}

    \item The three cases above cover all possibilities, so the verification is concluded.
\end{prooff}

Hence, the right $L(\R)$-approximations for super-suitable mice are suitable, short tree iterable premice.
The notion of a normal iterate for them is not in general defined, so we introduce a more general notion, which we shall call a \textit{normal pseudo-iterate}.

\begin{definition}\label{arrow}
    The relation \intro{$M\longrightarrow_\T^\mathrm{st}N$} holds iff
    \begin{parts}
        \item $M$ is suitable and short tree iterable,

        \item $\T$ is a tree on $M$ according to $\Sigma^\mathrm{st}_M$,

        \item\label{1109} if $\T$ has a last model, then $N$ is the last model of $\T$ and the $M$-to-$N$ branch of $\T$ does not drop,

        \item\label{1111} if $\T$ does not have a last model, then $\T$ is maximal and there exists a nondropping cofinal wellfounded branch $b$ through $\T$ such that $\M^\T_b$ is suitable and $\M^\T_b=N$.\qed
    \end{parts}
\end{definition}

This relation is absolute between $L(\R)$ and $V$.
In the part \ref{1109}, one has the behavior that one would expect from a normal iterate.
It is the part \ref{1111} that is unusual insomuch that we do not know what is the branch leading to the final model.
Let us look closer to this case.

\begin{lemma}\label{arrow and suitability}
    Suppose that $M\longrightarrow_\T^\mathrm{st}N$.
    Then $N$ is suitable and if $\T$ is maximal, then $\delta(\T)$ is the Woodin of $N$ and $N\para\delta(\T)=\M(\T)$.
\end{lemma}
\begin{proof}
    This is immediate from the definition since it subsumes the short tree iterability of $M$.
\end{proof}

What could happen is that there could be many branches that are ``good enough'', so we do not force ourselves to choose, but we know that they all lead to the same model.
In other words, by Lemma \ref{arrow and suitability}, we have that
$$N/\delta(\T)=\lp(\M(\T)),$$
which means that $N$ depends only on $\T$, but not on $b$.
That there are indeed such branches follows from the definition of the short tree iterability, but (assuming that $M$ is super-suitable) we do not know which one of them is chosen by $\Sigma_M$.
We observe below that in the case that $M$ is super-suitable, the relation
$$M\longrightarrow_\T^\mathrm{st}N$$
simply means that $N$ is a normal iterate of $M$ via $\T$.

\begin{lemma}\label{super-suitability and arrow}
    Suppose that
    \begin{assume}
        \item $M$ is super-suitable,

        \item $N$ is a countable premouse,

        \item $\T$ is a countable normal tree on $M$,

        \item\label{1146} for all limit $\gamma<\lh(\T)$, $\qq(\T\rest\gamma)$ exists.
    \end{assume}
    Then $M\longrightarrow_\T^\mathrm{st}N$ holds if and only if the following holds:
    \begin{parts}
        \item\label{1139} $\T$ is according to $\Sigma_M$,

        \item\label{1141} if $\T$ has a last model, then $\T\rest(\lh(\T)-1)$ is not maximal\footnote{We do not say ``short'' here because we $\T\rest(\lh(\T)-1)$ might not be of limit length.}, $N$ is the last model of $\T$, and the $M$-to-$N$ branch of $\T$ does not drop,

        \item\label{1143} if $\T$ does not have a last model, then for $b:=\Sigma_M(\T)$, we have that $b$ does not drop and $N=\M^\T_b$.
    \end{parts}
\end{lemma}
\begin{prooff}
    \item[] \textbf{Implication ($\Rightarrow$)}

    \item The part \ref{1139} follows from Proposition \ref{strategy and st strategy}, while the part \ref{1141} follows from the part \ref{1139}.

    \item We want to verify the part \ref{1143}.
    In this case, the tree $\T$ is maximal.

    \item If $b$ dropped, then $\qq(\T,b)$ would exist and be $(\omega_1+1)$-iterable, which would mean that $\qq(\T)$ exists and is equal to $\qq(\T,b)$.
    This being contradictory, we conclude that $b$ does not drop.

    \item It follows that $\M^\T_b$ is super-suitable and that $\delta(\T)$ is its Woodin.
    In particular,
    $$\M^\T_b/\delta(\T)=\lp(\M(\T))=N/\delta(\T),$$
    where the second equality follows from Lemma \ref{arrow and suitability}.

    \item This suffices for the conclusion that $N=\M^\T_b$.

    \item[] \textbf{Implication ($\Leftarrow$)}

    \item The assumptions yield that if $\T$ has a last model, then $\T\rest(\lh(\T)-1)$ is not maximal.

    \item We also know that there cannot exist $\xi$ such that $\xi+1<\lh(\T)$ and $\T\rest\xi$ is maximal.
    The reason for this is that $\delta(\T)=\delta(\M^\T_\xi)$ is a strong cutpoint of $\M^\T_\xi$, so any cofinal branch through $\T$ must contain $\xi$ and must drop at the successor of $\xi$.

    \item These two fact together with Proposition \ref{strategy and st strategy} imply that $\T$ is according to $\Sigma^\mathrm{st}_M$.
    The rest is now a routine verification.
\end{prooff}

\begin{definition}
    Suppose that $M$ is suitable and $N$ is a countable premouse.
    Then $N$ is a \intro{normal pseudo-iterate} of $M$ iff there exists $\T$ such that $M\longrightarrow_\T^\mathrm{st}N$.\qed
\end{definition}

\begin{proposition}\label{normal pseudo-iterate, normal iterate}
    Suppose that
    \begin{assume}
        \item $M$ is super-suitable,

        \item $N$ is a countable premouse.
    \end{assume}
    Then the following are equivalent.
    \begin{parts}
        \item $N$ is a normal pseudo-iterate of $M$.

        \item $N$ is a normal iterate of $M$.
    \end{parts}
\end{proposition}
\begin{proof}
    This follows from Lemma \ref{super-suitability and arrow}.
\end{proof}

It turns out that the short tree iterability suffices for the comparison.
In order to prove this, we require the following two lemmas.

\begin{lemma}\label{1269} \label{(12)l1}
    Assume $\mathrm{ZF}$.
    Then there exists a poset $\P \subseteq \omega_1$ such that $\P \in \hod$ and for every real $x$, there exists a $\hod$-generic filter $G$ for $\P$ such that
    $$\hod_x = \hod [G].$$
\end{lemma}
\begin{proof}
    This follows from \cite[Theorem 10.3]{larson2023extensions}.
\end{proof}

\begin{lemma}\label{1278} \label{(1)l2}
    Assume $\mathrm{ZF} + \mathrm{AD}^+$ and $V = L(\R)$.
    Then for every real $x$,
    $$\hod_x = \hod [x].$$
\end{lemma}
\begin{proof}
    See \cite[Corollary 10.24]{larson2023extensions}.
\end{proof}

\begin{proposition}\label{comparison} \label{(12)p}
    Suppose that $M,N$ are suitable and short tree iterable.
    Then there exists $P$ which is a normal pseudo-iterate of both $M$ and $N$.
\end{proposition}
\begin{prooff}
    \item Let $T$ be the tree of a $(\Sigma^2_1)^{L(\R)}$-scale on a universal $(\Sigma^2_1)^{L(\R)}$-set and let 
    $$\mu:=\max\{\delta(M),\delta(N)\}^{+L[T,M,N]}.$$
    We want to show that $\mu<\omega_1$, for once we do this, the conclusion follows from \cite[Theorem 6.14]{steel2016hod}.

    \item Let $x$ be a real coding $M$ and $N$.
    It suffices to show that $\omega_1^{L[T,x]}<\omega_1$.

    \item \label{1299} The theorem at the bottom of the page 77 of \cite{steel1995HOD} implies that 
    $$L[T]\cap V_{\boldsymbol\delta^2_1} = \hod \cap V_{\boldsymbol\delta^2_1}.$$
    We denote this set by $X$.
    Since $\boldsymbol\delta^2_1$ is an inaccessible cardinal in $\hod$, it holds that $X \models \mathrm{ZFC}$.

    \item By Lemmas \ref{1269} and \ref{1278}, $\hod_x = \hod [x]$, and this class is a ${<}\boldsymbol\delta^2_1$-generic extension of $\hod$.
    Consequently, $x$ is generic over $X$ and
    $$\hod_x \cap V_{\boldsymbol\delta^2_1} = \hod [x] \cap V_{\boldsymbol\delta^2_1} = \left( \hod \cap V_{\boldsymbol\delta^2_1} \right) [x] = X[x].$$

    \item By \ref{1299}, $L[T]\cap V_{\boldsymbol\delta^2_1} = X$ and $\boldsymbol\delta^2_1$ is a beth-fixed point in $L[T]$, so $x$ is ${<} \boldsymbol\delta^2_1$-generic over $L[T]$ and
    $$L[T][x] \cap V_{\boldsymbol\delta^2_1} = X[x].$$

    \item We can now conclude that
    $$\begin{aligned}
        \omega_1^{L[T, x]} &= \omega_1^{L[T][x]} = \omega_1^{X[x]} = \omega_1^{\hod_x} \\
        &< (\mbox{``the first measurable''})^{\hod_x} = \omega_1,
    \end{aligned}$$
    as required.
\end{prooff}

\section{HOD as a Direct Limit of Mice}
\label{HOD as a Direct Limit of Mice}

In this section, we describe $\hod\para\Theta$ as the direct limit of a certain directed family of mice.
This is just a slight reformulation of the well-known results of \cite{steel2016hod}.
In the end, we recall how these technics were used in \cite{neeman2007inner} to prove that there are in $L(\R)$ unique supercompactness measures on all $[\eta]^\omega$.
Here, we will go one step further than \cite{neeman2007inner} in that we will observe that these measures are, in fact, traces on $L(\R)$ of the club filters of $V$.

\begin{definition}
    Suppose that $M,N$ are super-suitable.
    Then we define the following.
    \begin{parts}
        \item \intro{$\ii(M)$} is the set of all normal iterates of $M$.

        \item \intro{$M\leq_\ii N$} iff $N$ is a normal iterate of $M$.\qed
    \end{parts}
\end{definition}

We note that if the main branch of a normal tree on a super-suitable mouse does not drop, then the tree is based below the Woodin.
The reason why extenders above the Woodin cannot be used is that once such an extender is used, all later extenders need to be above the Woodin, while any such extender necessarily leads to a drop.
Instead of restricting ourselves to normal iterates, we could have looked at all iterates (via stacks of normal trees).
However, by results of \cite{schlutzenberg2021full}, stacks of trees can be normalized, so nothing would be gained from this extension and we are still the scenario of \cite[Section 8]{steel2010outline} and \cite[Section 6]{steel2016hod}.
In particular, the structure $(\ii(M),\leq_\ii)$ is a countably directed partial order.

\begin{lemma}
    Suppose that $M$ is super-suitable.
    Then $(\ii(M),\leq_\ii)$ is a countably directed partial order.\qed
\end{lemma}

This partial order will index our directed family.
To every pair $(M,N)$ satisfying $M\leq_\ii N$, we need to assign an arrow $M\to N$, which is provided by the next lemma.

\begin{lemma}\label{(13)}
    Suppose that
    \begin{assume}
        \item $M$ is super-suitable,

        \item $N$ is a normal iterate of $M$.
    \end{assume}
    Then there exists a unique normal tree $\T$ on $M$ according to $\Sigma_M$ with the last model $N$.
    This tree is countable and its main branch does not drop.
\end{lemma}
\begin{proof}
    Such a tree $\T$ exists by the fact that $N$ is a normal iterate of $M$.
    The uniqueness follows from the fact that any such tree must be the first coordinate of the coiteration of $(M,\Sigma_M)$ and $(N,\Sigma_N)$.
    We will further elaborate on this last statement.
    A way to restate it is to say that $(\T, \uu)$ is the coiteration of $(M,\Sigma_M)$ and $(N,\Sigma_N)$ whenever $\T$ is a normal tree on $(M,\Sigma_M)$ whose last model is $N$ and $\uu$ is the trivial tree on $(N,\Sigma_N)$ whose length is $\lh(\T)$.\footnote{
    Recall that in coiterations we use padded trees, i.e. trees where one can pass from the $\alpha^\mathrm{th}$ model to the $(\alpha + 1)^\mathrm{th}$ model via the identity embedding.
    Having this in mind, a \textit{trivial tree} is a tree in which all embeddings are the identity (and the tree order is linear).
    }
    This is verified by induction.

    Let $\alpha < \lh(\T)$ be such that $(\T \rest \alpha, \uu \rest \alpha)$ is the coiteration up to stage $\alpha$.
    We want to show that $(\T \rest (\alpha + 1), \uu \rest (\alpha + 1))$ is the coiteration up to stage $\alpha + 1$.
    If $\alpha$ is $0$, this is trivially true; if $\alpha$ is a limit ordinal, this follows from the fact that $\T$ and $\uu$ are built using the same pair of stragies as in the case of the coiteration.
    So, let us assume that $\alpha$ is a successor ordinal.
    The least disagreement between $\M^\T_{\alpha - 1}$ and $\M^\uu_{\alpha - 1}$ is exactly the least disagreement between $\M^\T_{\alpha - 1}$ and $N$, which is just another way of saying  ``the least disagreement between $\M^\T_{\alpha - 1}$ and $\M^\T_{\lh(\T) - 1}$''.
    By \cite[Lemma 3.5]{steel2010outline}, this least disagreement occurs at the index $\lh(E^\T_{\alpha - 1})$, and this index is a cardinal in $\M^\T_{\lh(\T) - 1}$; in other words, the least disagreement is exactly $(E^\T_{\alpha - 1}, \emptyset)$.\footnote{
    To understand the second coordinate, recall that no extender may be indexed at a cardinal.
    }
    This means that to get from $(\T \rest \alpha, \uu \rest \alpha)$ to $(\T \rest (\alpha + 1), \uu \rest (\alpha + 1))$, one indeed uses the least disagreement between $\M^\T_{\alpha - 1}$ and $\M^\uu_{\alpha - 1}$, so $(\T \rest (\alpha + 1), \uu \rest (\alpha + 1))$ remains an initial segment of the coiteration.

    The argument just given establishes that $(\T, \uu)$ is an initial segment of the coiteration of $(M,\Sigma_M)$ and $(N,\Sigma_N)$.
    However, the final models of $\T$ and $\uu$ are equal (viz., $N$), so the coiteration terminates exactly at this stage.
    This shows that $(\T, \uu)$ is exactly the coiteration of $(M,\Sigma_M)$ and $(N,\Sigma_N)$.
\end{proof}





\begin{definition}
    Suppose that
    \begin{assume}
        \item $M$ is super-suitable,

        \item $N$ is a normal iterate of $M$.
    \end{assume}
    Then the mapping \intro{$\pi_{M,N}$} is defined to be the mapping 
    $$\pi^\T:M\to N$$
    where $\T$ is the normal tree on $M$ according to $\Sigma_M$ whose last model is $N$.\qed
\end{definition}

We are ready to introduce the directed system.

\begin{definition}
    Suppose that $M$ is super-suitable.
    Then we define \intro{$\Vec\ii(M)$} to be the system
    $$(P,\pi_{P,Q} : P,Q\in\ii(M),\ P\leq_\ii Q).$$
    \qed
\end{definition}

\begin{proposition}
    Suppose that $M$ is super-suitable.
    Then $\Vec{\ii}(M)$ is a directed system with a wellfounded direct limit.
\end{proposition}
\begin{proof}
    The commutativity follows from the normalization (cf. \cite[Theorem 1.1]{schlutzenberg2021full}), while wellfoundedness follows from the fact that $(\ii(M),\leq_\ii)$ is countably directed.
\end{proof}

By the results of Section \ref{section: short tree iterability}, and in particular Proposition \ref{comparison}, if $M$ and $N$ are two super-suitable mice, then there exists a super-suitable $P$ which is a normal iterate of both.
We get that
$$\Vec{\ii}(M)\rest P=\Vec{\ii}(N)\rest P=\Vec{\ii}(P).$$
This leads to the following consequence.

\begin{proposition}
    Suppose that $M,N$ are super-suitable.
    Then the direct limits of the systems $\Vec{\ii}(M)$ and $\Vec{\ii}(N)$ are equal and for all $P\in\ii(M)\cap\ii(N)$,
    $$\pi^{\Vec{\ii}(M)}_{P,\infty}=\pi^{\Vec{\ii}(N)}_{P,\infty}.$$
    \qed
\end{proposition}

\begin{definition}\label{1630'}
    Let us denote $\M:=\M_\omega\parallel\delta(\M_\omega)^{+\M_\omega}$.
    We define objects \intro{$\M_\infty$, $\pi_{P,\infty}$} (for all $P\in\ii(\M)$) as follows:
    $$(\M_\infty,\pi_{P,\infty} : P\in\ii(\M))$$
    is the direct limit of the system $\Vec{\ii}(\M)$.\qed
\end{definition}

We can now state the characterization of $\hod\para\Theta$ in the style of \cite{steel2016hod} that we shall use.

\begin{theorem}\label{1639}
    It holds that $\delta(\M_\infty)=\Theta$ and $\M_\infty\para\Theta=\hod\para\Theta$.\qed
\end{theorem}

The paper \cite{neeman2007inner} used this representation of $\hod \para \Theta$ to define certain filters on $[\eta]^\omega$ (all $\eta \in [\omega, \Theta)$) in $L(\R)$ which turn out to be the unique supercompactness measures on these sets.
We will be able to identify these measure as the traces on $L(\R)$ of the club filters from $V$.

\begin{lemma}\label{1889}
    Let $\eta < \Theta$.
    Then $\eta^\omega \subseteq L(\R)$.
\end{lemma}
\begin{proof}
    Let $\sigma \in \eta^\omega$ be given.
    Fix a set of reals $A \in L(\R)$ of Wadge rank $\eta$.
    For every $n < \omega$, there is a real $x_n$ which codes a continuous reduction of the set $A$ to a set of reals of Wadge rank $\sigma_n$.
    Let $y$ be a real that codes the sequence $(x_n : n < \omega)$.
    Now, the pair $(A, y)$ belongs to $L(\R)$ and codes the sequence $\sigma$ in the following way: from $y$, one decodes the sequence of reals $(x_n : n < \omega)$, then the sequence of continuous functions $(f_n : n < \omega)$ from $\R$ to $\R$, and finally, the sequence
    $$\sigma = (f_n^{-1}[A] : n < \omega).$$
    Therefore, $\sigma \in L(\R)$.
\end{proof}

The lemma is relevant for us because we get $[\eta]^\omega \in L(\R)$.
Looking now at any filter on $[\eta]^\omega$ in $V$, its intersection with $L(\R)$ becomes an $L(\R)$--filter on $[\eta]^\omega$.

\begin{definition}
    Let $\eta \in [\omega, \Theta)$.
    Then \intro{$\mu_\eta$} denotes the filter $\mathrm{club}_{[\eta]^\omega} \cap L(\R)$, where $\mathrm{club}_{[\eta]^\omega}$ is the club filter on $[\eta]^\omega$.
    \qed
\end{definition}

\begin{theorem}\label{1896}
    Let $\eta \in [\omega, \Theta)$.
    Then $\mu_\eta$ belongs to $L(\R)$ and it is the unique supercompactness measure\footnote{
    We remind the reader that a supercompactness measure on $[X]^\omega$ is just a normal fine ultrafilter on that set.
    } on $[\eta]^\omega$ in that model.
\end{theorem}
\begin{prooff}
    \item By \cite[Corollary 4.14]{neeman2007inner}, there exists in $L(\R)$ a unique supercompactness measure $\mu$ on $[\eta]^\omega$.
    We will be done if we show that $\mu = \mu_\eta$.

    \item Since $\mu$ is an $L(\R)$-ultrafilter, it suffices to show that every $S \in \mu$ contains a club (in $V$).

    \item For $x \in \R$, we denote by $\M_x$ the $x$-premouse $\M_\omega(x) \para \tau$, where $\tau := \delta(\M_\omega(x))^{+\M_\omega(x)}$, and by $\Sigma_{\M_x}$ the $(\omega_1, \omega_1 + 1)$-iteration strategy for $\M_x$ derived from a fixed $(\omega_1, \omega_1 + 1)$-iteration strategy for $\M_\omega(x)$.
    An \textit{$x$-super-suitable premouse} is a $\Sigma_{\M_x}$-iterate of $\M_x$.
    For an $x$-super-suitable premouse $M$, we denote by $\Sigma_M$ the strategy for $M$ induced by $\Sigma_{\M_x}$.

    \item By repeating the definitions given above while replacing $\M$ and $\Sigma_\M$ by $\M_x$ and $\Sigma_{\M_x}$ (resp.), super-suitable by $x$-super-suitable, etc., we obtain the following objects:
    \begin{itemize}
        \item $(\ii_x(M), \leq_{\ii_x})$,

        \item $\pi^x_{M, N}: M \to N$,

        \item $\Vec{\ii}_x(M)$,

        \item $\M_\infty^x$,

        \item $\pi^x_{M, \infty}$,
    \end{itemize}
    where $x$ is a real, $M$ is $x$-super-suitable, and $N$ is a normal iterate of $M$.
    Furthermore, the facts of this section remain true when conditioned on a real $x$ (and when $\hod$ is replaced by $\hod_x$).



    \item For $x \in \R$, a sequence of the form $(M_n, \T_n: n < \omega)$ is a \textit{$x$-super-correct iteration sequence} iff $M_0$ is $x$-super-suitable and for all $n < \omega$, $\T_n$ is a countable normal tree on $M_n$ according to $\T_n$ with the last model $M_{n + 1}$; we say that this $x$-super-correct iteration sequence is \textit{on $M_0$}.
    For an $x$-super-correct iteration sequence $\sigma$, we denote by $a(\sigma)$ the countable set
    $$\bigcup_{n < \omega} \ran (\pi^x_{M_n, \infty}) \cap \eta.$$

    \item For $x \in \R$ and for an $x$-super-suitable $M$, let $C^x_M$ consist of all $a(\sigma)$ where $\sigma$ is an $x$-super-correct iteration sequence on $M$.
    By \cite[Remark 4.11]{neeman2007inner}, there exist $x \in \R$ and an $x$-super-suitable $M$ such that $C^x_M \subseteq S$.

    \item We will be done if we show that $C^x_M$ contains a club.
    By \cite[Proposition 3.14]{kasum2024building}, it suffices to show that Player II has a winning quasi-strategy in $\Game^{\mathrm{club}}(\eta, C^x_M)$.

    \item Let $\Game^*$ be a length $\omega$ two player game of the form
    \begin{center}
        \begin{tabular}{c|ccccc}
            I  & $\alpha_0$ &                     & $\alpha_1$ &                          & $\cdots$\\
            \hline
            II &       & $\beta_0$, $\T_0$, $M_1$ &            & $\beta_1$, $\T_1$, $M_2$ & $\cdots$
        \end{tabular}
    \end{center}
    where $M_0 := M$ by convention and for all $n < \omega$,
    \begin{itemize}
        \item $\alpha_n, \beta_n < \eta$,

        \item $\T_n$ is a normal tree on $M_n$ according to $\Sigma_{M_n}$ with the last model $M_{n + 1}$.
    \end{itemize}
    Player II wins iff
    $$a \left( (M_n, \T_n: n < \omega) \right) = \{\alpha_n, \beta_n: n < \omega\}.$$

    \item Observe that for every play $(\alpha_n, (\beta_n, \T_n, M_{n + 1}): n < \omega)$ of $\Game^*$ which is won by Player II, $(\alpha_n, \beta_n: n < \omega)$ is a play of $\Game^\mathrm{club}(\eta, C^x_M)$ which is won by Player II.
    So, it suffices to show that Player II has a winning quasi-strategy in the game $\Game^*$.

    \item This readily follows: a winning quasi-strategy for Player II in $\Game^*$ consists in
    \begin{itemize}
        \item ensuring that for all $n < \omega$, $\alpha_n \in \ran(\pi_{M_{n + 1}, \infty})$,

        \item fixing an enumeration $\{\gamma^n_i: i < \omega\}$ of $\ran(\pi^x_{M_n, \infty}) \cap \eta$ for every $n < \omega$,

        \item using a bookkeeping device when choosing $\beta_n$'s so as to ensure that
        $$\{\beta_n : n < \omega\} = \{\gamma^n_i: n, i < \omega\}.$$
    \end{itemize}
    Player II indeed wins in this way: in the end, we have
    $$\{\alpha_n, \beta_n : n < \omega\} = \{\gamma^n_i: n, i < \omega\} = a \left( (M_n, \T_n : n < \omega) \right).$$
\end{prooff}


\section{\texorpdfstring{$\hod$}{HOD} as a Normal Iterate}
\label{2031}

Consider a super-suitable mouse $M$.
The results of \cite{schlutzenberg2021full} imply that there is a non-dropping normal tree on $M$ whose last model is $\M_\infty$.
We will argue in this section for the uniqueness of this tree, provided it satisfies the additional requirement that its final branch embedding is equal to the direct limit embedding $\pi_{M, \infty}$.
It is the study of the properties of this tree that will eventually lead to the local definability of $\hod$ in $L(\R)$.

The requirement that the tree embedding is equal to $\pi_{M, \infty}$ will only be used to argue that the last branch must be chosen in a unique way.
The earlier branches are chosen according to Q-structures.
The next lemma establishes existence of these Q-structures.

\begin{lemma}\label{1997}
    Let $M$ be super-suitable and, in a generic extension of $V$, let $\T$ be a normal tree on $M$ whose last model is $\M_\infty$.\footnote{
    Here and in similar statements below, we understand that $\M_\infty$ means $\M_\infty^V$.
    The point is that we are considering the trees that end with $\M_\infty^V$, but may belong to some generic extension of $V$.
    }
    Then for all $\alpha < \lh(\T) - 1$, either $\M^\T_\alpha$ does not have a Woodin cardinal or $\lh(E^\T_\alpha)$ is strictly less than the least Woodin cardinal of $\M^\T_\alpha$.
\end{lemma}
\begin{prooff}
    \item Recall that $M$ is an initial segment of an iterate of $\M_\omega^\sharp$.
    This means that no extender on the $M$ sequence, or on the sequence of any of its iterates, contains an extender overlapping a Woodin cardinal.

    \item Assume towards a contradiction that the conclusion fails and let $\alpha < \lh(\T) - 1$ be the least such that $\delta(\M^\T_\alpha)$ is defined and $\lh(E^\T_\alpha) \geq \delta(\M^\T_\alpha)$.

    \item By \cite[Lemma 3.5]{steel2010outline}, we get that $\M^\T_\alpha$ and the last model of $\T$ (viz., $\M_\infty$) agree below $\lh(E^\T_\alpha)$ and $\lh(E^\T_\alpha)$ is a cardinal in $\M_\infty$.
    
    \item Since $E^\T_\alpha$ does not overlap $\delta(\M^\T_\alpha)$, we have that $\crit(E^\T_\alpha) \geq \delta(\M^\T_\alpha)$.

    \item The last two points yield that in $\M_\infty$, there is a measurable Woodin cardinal or there is a measurable cardinal above a Woodin cardinal -- either case is absurd.
\end{prooff}





Since the intermediate branches are chosen according to Q-structures, we must see how to identify them.
We will read them off of $\M_\infty$.
The next lemma helps us understand that the relevant part of $\M_\infty$ is the part up to $\Theta$.

\begin{lemma}\label{2095}
    Let $M$ be super-suitable and, in a generic extension of $V$, let $\T$ be a normal tree on $M$ whose last model is $\M_\infty$.
    Then for all $\xi$ such that $\xi + 1 < \lh(\T)$, it holds that $\lh(E^\T_\xi) < \Theta$.
\end{lemma}
\begin{proof}
    Assume towards contradiction that $\lh(E^\T_{\xi}) \geq \Theta$.
    By \cite[Lemma 3.5]{steel2010outline}, $\lh(E^\T_\xi)$ is a successor cardinal in the last model of $\T$, i.e. in $\M_\infty$.
    This is a contradiction: $\M_\infty$ has no successor cardinals which are $\geq \Theta$.
\end{proof}

Thus, if the information contained $\M_\infty$ at level $\Theta$ or higher is relevant, then we are at the last stage of the iteration from $M$ to $\M_\infty$.
In this last stage, the height of the lined-up part is equal to $\Theta$.

\begin{corollary}\label{2030}
    Let $M$ be super-suitable and, in a generic extension of $V$, let $\T$ be a normal tree on $M$ whose last model is $\M_\infty$ and whose main branch does not drop.
    Suppose $\lambda < \lh(\T)$ is a limit ordinal such that $\delta(\T \rest \lambda) \geq \Theta$.
    Then $\lh(\T) = \lambda + 1$ and $\delta(\T \rest \lambda) = \Theta$.
\end{corollary}
\begin{proof}
    By Lemma \ref{2095}, $\delta(\T \rest \lambda) \leq \Theta$, so in fact, $\delta(\T \rest \lambda) = \Theta$.
    Any extender after stage $\lambda$ would have the index $\geq \Theta$, which is excluded as a possibility by Lemma \ref{2095}.
    Therefore, $\lh(\T) = \lambda + 1$.
\end{proof}

We can now conclude that the intermediate branches of the iteration to $\M_\infty$ are determined by Q-structures.

\begin{lemma}\label{2030'}
    Let $M$ be super-suitable and, in a generic extension of $V$, let $\T$ be a normal tree on $M$ whose last model is $\M_\infty$.
    Suppose $\lambda < \lh(\T)$ is a limit ordinal such that $\delta(\T \rest \lambda) < \Theta$.
    Then $[0, \lambda)_{\T}$ is the unique well-founded cofinal branch $b$ through $\T \rest \lambda$ such that $\qq(\T \rest \lambda, b)$ exists and is an initial segment of $\hod \para \Theta$.
\end{lemma}
\begin{prooff}
    \item We only need to show that $\qq(\T \rest \lambda, b)$ exists and is an initial segment of $\hod \para \Theta$.
    The uniqueness follows from Proposition \ref{252}.
    
    \item Let $\delta:=\delta(\T \rest \lambda)=\bigcup_{i<\lambda}\lh(E^{\T}_i)$.
    So, we are assuming that $\delta<\Theta$.

    \item If $[0, \lambda)_{\T}$ drops, then $\qq(\T \rest \lambda, [0, \lambda)_{\T})$ exists.

    \item \claim Even if $[0, \lambda)_{\T}$ does not drop, $\qq(\T \rest \lambda, [0, \lambda)_{\T})$ exists.
    \begin{proof}
        It suffices to prove that $\delta$ strictly less than the least Woodin cardinal of $\M^\T_\lambda$, viz. $\delta(\M^\T_\lambda)$.
        Indeed, then $\delta$ is not a Woodin cardinal in $\M^\T_\lambda$, so $\qq(\T \rest \lambda, [0, \lambda)_{\T})$ exists.

        Notice that
        $$\delta = \bigcup_{\alpha \in [0, \lambda)_\T} \lh(E^\T_\alpha).$$
        Now, by Lemma \ref{1997}, we get that for all $\alpha \in [0, \lambda)_\T$,
        $$\lh(E^\T_\alpha) < \delta(\M^\T_\alpha) \leq i^\T_{\alpha, \lambda}(\delta(\M^\T_\alpha)) = \delta(\M^\T_\lambda),$$
        so we conclude that $\delta \leq \delta(\M^\T_\lambda)$.

        Assume towards contradiction that $\delta = \delta(\M^\T_\lambda)$.
        Then 
        $$\delta(\M^\T_\lambda) < \Theta = \delta(\M_\infty),$$
        which means that $\M^\T_\lambda$ is not the last model of $\T$.
        Hence, there is $E^\T_\lambda$ and by Lemma \ref{1997}:
        $$\lh(E^\T_\lambda) < \delta (\M^\T_\lambda) = \delta = \delta(\T \rest \lambda).$$
        This is, of course, a contradiction.

        Therefore, $\delta < \delta(\M^\T_\lambda)$, as required.
    \end{proof}

    \item Thus, in any case, $\qq(\T \rest \lambda,[0, \lambda)_{\T})$ exists.
    Let us denote it by $Q$.
    Our goal is to show that $Q\unlhd\hod\para\Theta$.

    \item Since $\T$ is normal, $\lh(E^{\T}_\lambda) \geq \delta$.

    \item Since for all $i<\lambda$, $\lh(E^{\T}_i)$ is a cardinal of $\M^{\T}_\lambda$, we know that that $\delta$ is a limit cardinal of $\M^{\T}_\lambda$.
    This means that there are no extenders in $\M^{\T}_\lambda$ indexed at $\delta$ and consequently, it must be the case that $\lh(E^{\T}_\lambda)$ is strictly greater than $\delta$.

    \item \label{x(16)l1} Since $\lh(E^{\T}_\lambda)$ is a cardinal in $\hod$,\footnote{
    See \cite[Lemma 3.5]{steel2010outline}.
    } we conclude that
    $$\delta^{+\hod} \leq \lh(E^{\T}_\lambda).$$

    \item \label{x1721} Since $\M^{\T}_\lambda \para \lh(E^{\T}_\lambda) = \hod \para \lh(E^{\T}_\lambda)$,\footnote{
    This is again follows from \cite[Lemma 3.5]{steel2010outline}.
    } it follows that
    $$\hod \para \delta^{+\hod} \unlhd \M^{\T}_\lambda.$$


    \item Since $\delta$ is not Woodin\footnote{This is because $\delta<\Theta$, the models $\hod$ and $\M_\infty$ agree below $\Theta$, and $\Theta$ is the only Woodin of $\M_\infty$.} in $\hod$, there exists the smallest $\zeta$ such that 
    $$\hod\para(\zeta+1)\models``\delta\mbox{ is not Woodin}".$$

    \item \label{x(16)l2} Observe that $\zeta<\delta^{+\hod}$, so \ref{x1721} implies that 
    $$\hod\para(\zeta+1) = \M^{\T}_\lambda \para (\zeta + 1).$$
    Therefore, $Q=\hod\para\zeta$, as required.
\end{prooff}

Let us reformulate the lemma more clearly.

\begin{proposition}\label{2013}
    Suppose that $M$ is super-suitable.
    Then in every generic extension of $V$, for all normal trees $\T$ on $M$ whose last model is $\M_\infty$, it holds that:
    \begin{parts}
        \item \label{2017} for all $i < \lh(\T) - 1$, $\lh(E^\T_i) < \Theta$ and is equal to the least $\alpha < \delta(\M^\T_i)$ such that $\M^\T_i | \alpha \neq \hod | \alpha$;

        \item \label{2019} for all limit $\lambda < \lh(\T) - 1$, $[0, \lambda)_{\T}$ is the unique well-founded cofinal branch $b$ through $\T \rest \lambda$ such that $\qq(\T \rest \lambda, b)$ exists and is an initial segment of $\hod \para \Theta$.
    \end{parts}
\end{proposition}
\begin{proof}
    Part \ref{2017} is immediate once we have Lemmas \ref{1997} and \ref{2095}.
    On the other hand, part \ref{2019} is verified as follows: Corollary \ref{2030} yields that $\delta(\T \rest \lambda) < \Theta$; the rest follows from Lemma \ref{2030'}.
\end{proof}

In the next lemma we finish verification that there is at most one iteration to $\M_\infty$.

\begin{lemma}\label{2024}
    Suppose that $M$ is super-suitable.
    Then in every generic extension of $V$, there is at most one normal tree on $M$ whose last model is $\M_\infty$, whose main branch does not drop and whose main-branch embedding is equal to $\pi_{M, \infty}$.
\end{lemma}
\begin{prooff}
    \item Let us work in a fixed generic extension of $V$ and let $\T$ and $\uu$ be trees on $M$ whose last models are both $\M_\infty$, whose main branches do not drop and whose main-branch embeddings are equal to $\pi_{M, \infty}$.
    We want to show that $\T = \uu$.

    \item We will prove by induction on $\alpha \leq \min\{\lh(\T), \lh(\uu)\}$ that $\T \rest \alpha = \uu \rest \alpha$.
    Once we do this, we will know that either $\T \unlhd \uu$ or $\uu \unlhd \T$.
    Since $\T$ and $\uu$ end with the same model, it follows that $\T = \uu$.

    \item Let us proceed with the induction.
    We are assuming that $\alpha \leq \min\{\lh(\T), \lh(\uu)\}$ and that for all $\beta < \alpha$, $\T \rest \beta = \uu \rest \beta$ and we need to show that $\T \rest \alpha = \uu \rest \alpha$.

    \item If $\alpha$ is $0$ or $1$ or a limit ordinal, there is nothing to be proved.
    Hence, the remaining cases are ``$\alpha$ is the successor of a limit ordinal'' and ``$\alpha$ is the successor of a successor ordinal''.
    In either case, let us denote $\beta := \alpha - 1$.

    \item \textbf{Case 1.} $\alpha$ is the successor of a limit ordinal.
    \begin{proof}
        In this case, we have to show that $[0, \beta)_\T = [0, \beta)_\uu$.
        If $\delta(\T \rest \beta) = \delta(\uu \rest \beta) < \Theta$, then by Lemma \ref{2030'}, both $[0, \beta)_\T$ and $[0, \beta)_\uu$ are the unique well-founded cofinal branch $b$ through the tree $\T \rest \beta = \uu \rest \beta$ which has a Q-structure which is an initial segment of $\hod \para \Theta$; in particular, $[0, \beta)_\T = [0, \beta)_\uu$.
        
        It remains to verify the conclusion in the case that the ordinal $\delta(\T \rest \beta) = \delta(\uu \rest \beta) = \Theta$.
        By Corollary \ref{2030}, $\beta$ is the last index in $\T$ and $\uu$, i.e.
        $$\lh(\T) = \lh(\uu) = \beta + 1 = \alpha.$$
        This means that the branches $[0, \beta)_\T$ and $[0, \beta)_\uu$ do not drop and $i^\T_{0, \beta} = \pi_{M, \infty}$ and $i^\uu_{0, \beta} = \pi_{M, \infty}$; in other words, branches $[0, \beta)_\T$ and $[0, \beta)_\uu$ through the tree $\T \rest \beta = \uu \rest \beta$ give rise to the same embedding.
        \cite[Lemma 2.6]{steel2016hod} now implies that these two branches are equal.
    \end{proof}

    \item \textbf{Case 2.} $\alpha$ is the successor of a successor ordinal.
    \begin{proof}
        In this case, both $\T \rest \alpha$ and $\uu \rest \alpha$ are obtained from the tree $\T \rest \beta = \uu \rest \beta$ by applying the same extender (Proposition \ref{2013}).
        Therefore, $\T \rest \alpha = \uu \rest \alpha$.
    \end{proof}

    \item The two cases conclude the induction.
\end{prooff}

It remains to see that there is at least one iteration to $\M_\infty$.
This was done in \cite{schlutzenberg2021full}.
We put everything together in the next proposition.

\begin{proposition}\label{2033}
    Suppose that $M$ is super-suitable.
    Then there exists a unique normal tree $\T$ on $M$ whose last model is $\M_\infty$, whose main branch does not drop and whose main-branch embedding is equal to $\pi_{M, \infty}$.
\end{proposition}
\begin{prooff}
    \item Taking into account Lemma \ref{2024}, we only need to verify existence.
    
    \item Let $\Psi$ be the restriction of $\Sigma_M$ to the stacks of normal trees in $H_{\omega_1}$ which have no unnecessary drops.
    This $(\omega_1, \omega_1)^*$-iteration strategy for $M$ is coded by a universally Baire by Theorem \ref{47}.
    
    \item Let $g$ be $V$-generic for $\col(\omega, \R)$ and let us work in $V[g]$.
    We have that $\Psi_g \cap H_{\omega_1}^V = \Psi$.\footnote{
    The strategy $\Psi$ is coded by a ${<}\lambda$-universally Baire set of reals, where $\lambda$ is the least limit of Woodin cardinals (see Theorem \ref{47}).
    This set of reals has the canonical reinterpretation in $V[g]$ and this reinterpretation defines an $(\omega_1, \omega_1)$-iteration strategy in $V[g]$, which we denote by $\Psi_g$.
    }

    \item \claim There exists a stack $(\T_i)_{i < \omega}$ on $(\M, \Psi_g)$ whose last model is $\M_\infty$, whose cofinal branch does not drop and gives rise to the embedding $\pi_{M, \infty}$.
    \begin{proof}
        In $V$, $\R$ surjects onto the directed system for $\M_\infty$.
        Hence, this directed system is countable in $V[g]$.
        It follows that there is a cofinal sequence of indices $(M_i : i < \omega)$ in the directed system, satisfying $M_0 = M$.
        In $V$, each $M_{i + 1}$ is the last model of a normal tree $\T_i$ on $M_i$ according to $\Sigma_{M_i}$; now, $\Sigma_{M_i}$ is a tail strategy of $\Sigma_M$, so $\T_i$ is according to $\Psi$.
        We have thus obtained a stack $(\T_i: i < \omega)$ on $M$ according to $\Psi_g$ with all desired properties.
    \end{proof}

    \item \claim There exists an $(\omega_1, \omega_1 + 1)^*$-iteration strategy $\mathrm{X}$ for $M$ such that $\Psi_g \subseteq \mathrm{X}$.
    \begin{proof}
        Let $s$ be a stack on $(M, \Psi_g)$ of total length $\omega_1$.
        We want to define $\mathrm{X}(s)$.
        Let $h$ be $V[g]$-generic for $\col(\omega, \omega_1)$.
        The stack $s$ is on $(M, \Psi_{g*h})$, so $b := \Psi_{g*h}(s)$ is a cofinal well-founded branch through $s$.
        Since $\col(\omega, \omega_1)$ is homogeneous in $V[g]$, we have that $b \in V[g]$.
        We may define $\mathrm{X}(s) := b$.
    \end{proof}

    \item By \cite[Theorem 1.1]{schlutzenberg2021full} applied to $(M, \mathrm{X})$, there is a normal tree $\T$ on $(M, \Psi_g)$ whose last model is $\M_\infty$, whose main branch does not drop, and whose main-branch embedding is equal to that of the stack $(\T_i)_{i < \omega}$, i.e. to $\pi_{\M, \infty}$.

    \item By Lemma \ref{2024}, $\T$ is the unique normal tree on $M$ whose last model is $\M_\infty$, whose main branch does not drop and whose main-branch embedding is equal to $\pi_{M, \infty}$.

    \item Homogeneity of $\col(\omega, \R)$ implies that $\T \in V$, which concludes the proof.
\end{prooff}

\begin{definition}
    Suppose that $M$ is super-suitable.
    Then we define \intro{$\T_M$} to be the unique normal tree on $M$ such that $\M_\infty$ is its last model, its main branch does not drop and $\pi^\T = \pi_{M, \infty}$.
    \qed
\end{definition}

\section{\texorpdfstring{$\hod$ as a Normal Iterate in $L(\R)$}{HOD as a Normal Iterate in L(R)}}

For every super-suitable mouse $M$, we have defined a normal iteration $\T_M$ on $M$.
We now want to establish that the length of this iteration is $\Theta + 1$, that its part up to $\Theta$ belongs to $L(\R)$, and that its final branch does \textbf{not} belong to $L(\R)$.

We start by establish an upper bound on the length of $\T_M$.

\begin{lemma}\label{2044}
    Let $M$ be super-suitable.
    Then $\lh(\T_M) \leq \Theta+1 < \cont^+$.
\end{lemma}
\begin{proof}
    By Lemma \ref{2095}, for all $i < \lh(\T_M) - 1$, we have that $\lh(E^{\T_M}_i) < \Theta$; in other words, we have the mapping
    $$\lh(\T_M) - 1 \to \Theta: i \mapsto \lh(E^{\T_M}_i).$$
    Since $\T_M$ is normal, this mapping is strictly increasing.
    Therefore, $\lh(\T_M) - 1 \leq \Theta$.
    To finish, notice that $\M_\infty$ is a surjective image of $\R$ (which easily from its definition as a direct limit).
    Since $\Theta < o(\M_\infty)$, we conclude that $|\Theta| \leq \cont$.
\end{proof}

The next thing to verify is that for all $\eta < \lh(\T_M)$, the tree $\T_M \rest \eta$ belongs to $L(\R)$.
This is done by induction.
The only non-trivial part is to show that if all proper initial segments of $\T_M \rest \eta$ belong to $L(\R)$, then so does $\T_M \rest \eta$.
In this step, the point is to characterize the branch $[0, \eta)_{\T_M}$ in a way which will make it easy to see that it belongs to $L(\R)$.

The first step towards this goal is to show that the intermediate branches of $\T_M$ are chosen according to \textbf{countably iterable} Q-structures.
We will show more: these Q-structures are actually $(\cont^+ + 1)$-iterable.
This is because they are initial segments of $\M_\infty$, for which we can find iteration strategies for large trees thanks to the universally Baire strategy for $\M_\omega^\sharp$.

\begin{lemma}\label{2033'}
    Let $\kappa$ be the least inaccessible cardinal.
    Then $\M_\omega^\sharp$ is fully iterable in $V_\kappa$.
\end{lemma}
\begin{proof}
    By Theorem \ref{47}, $\M_\omega^\sharp$ has a unique $(\omega_1, \omega_1 + 1)$-iteration strategy $\Psi$, and $\Psi \rest H_{\omega_1}$ is coded by a ${<}\kappa$-universally Baire set.
    This means that the following function $\Psi^+$ is a $\kappa$-iteration strategy for $\M_\omega^\sharp$ that is definable over $V_\kappa$: if $\mathcal{X}$ is a normal tree on $\M_\omega^\sharp$ of limit length according to $\Psi^+$, then $\Psi^+(\mathcal{X})$ is the unique $b$ such that
    $$V^{\col(\omega, \lh(\mathcal{X}))} \models \Psi_{\dot g}(\check{\mathcal{X}}) = \check{b};$$
    here, $\Psi_{\dot g}$ is the strategy coded by the canonical reinterpretation in $V^{\col(\omega, \lh(\mathcal{X}))}$ of the ${<}\kappa$-universally Baire set coding $\Psi$.
\end{proof}

\begin{lemma}\label{2044'}
    Let $\kappa$ be the least inaccessible cardinal.
    Then:
    \begin{parts}
        \item there exists a unique $\kappa$-iteration strategy $\Sigma$ for $\M_\omega \para \kappa$ which is definable over $V_\kappa$;

        \item letting $\M_\infty^+$ be the direct limit of all normal iterates of $(\M_\omega \para \kappa, \Sigma)$ which are based on the first Woodin of $\M_\omega \para \kappa$, we have that $\M_\infty$ is the initial segment of $\M_\infty^+$ cut at the successor of its first Woodin cardinal;

        \item letting $\pi: \M_\omega \para \kappa \to \M_\infty^+$ be the direct limit embedding, we have that there is a non-dropping normal tree $\T$ on $\M_\omega \para \kappa$ with the last model $\M_\infty^+$ which is according to $\Sigma$ and satisfies $\pi^\T = \pi$.
    \end{parts}
\end{lemma}
\begin{proof}
    Thanks to Lemma \ref{2033'}, can apply \cite[Corollary 1.2]{schlutzenberg2021full} to $V_\kappa$.
    This application yields that there is a unique $\kappa$-iteration strategy $\Sigma$ for $\M_\omega \para \kappa$ which is definable over $V_\kappa$\footnote{
    The point here is that we are the stated corollary inside $V_\kappa$, so ``definable'' means that it is a class in that universe.
    } and that there is a non-dropping normal tree $\T$ on $\M_\omega \para \kappa$ with the last model $\M_\infty^+$ which is according to $\Sigma$ and satisfies $\pi^\T = \pi$.\footnote{\cite[Corollary 1.2]{schlutzenberg2021full} does not state the latter part explicitly, but it does hold because the tree in question comes from \cite[Theorem 1.1]{schlutzenberg2021full} (where meeting of this requirement is stated).}
    It is immediate that $\M_\infty$ is the initial segment of $\M_\infty^+$ cut at the successor of its first Woodin cardinal.
\end{proof}

\begin{proposition}\label{2011}
    $\M_\infty$ is $(\cont^+ + 1)$-iterable.
\end{proposition}
\begin{prooff}
    \item Let $\kappa$ be the least inaccessible cardinal.
    Lemma \ref{2044'} implies
    \begin{parts}
        \item there exists a unique $\kappa$-iteration strategy $\Sigma$ for $\M_\omega \para \kappa$ which is definable over $V_\kappa$;

        \item letting $\M_\infty^+$ be the direct limit of all normal iterates of $(\M_\omega \para \kappa, \Sigma)$ which are based on the first Woodin of $\M_\omega \para \kappa$, we have that $\M_\infty$ is the initial segment of $\M_\infty^+$ cut at the successor of its first Woodin cardinal;

        \item letting $\pi: \M_\omega \para \kappa \to \M_\infty^+$ be the direct limit embedding, we have that there is a non-dropping normal tree $\T$ on $\M_\omega \para \kappa$ with the last model $\M_\infty^+$ which is according to $\Sigma$ and satisfies $\pi^\T = \pi$.
    \end{parts}
    




    \item Hence, we will be done if we show that $\M_\infty^+$ is $(\cont^+ + 1)$-iterable.
    To show it, we will exhibit a partial $(\cont^+, \cont^+ + 1)$-iteration strategy $\Tau$ for $\M^+$ such that:
    \begin{parts}
        \item $\Tau$ acts on the stacks with no artificial drops;

        \item $\T$ is according to $\Tau$.
    \end{parts}
    This will imply that $\M_\infty^+$ is $(\cont^+ + 1)$-iterable because the fragment of $\Tau_\T$ acting on normal trees is a $(\cont^+ + 1)$-iteration strategy for it.

    \item By \cite[Theorem 1.1]{schlutzenberg2021full}, there exists a partial $(\cont^+, \cont^+ + 1)$-iteration strategy $\Tau$ for $\M^+$ such that:
    \begin{parts}
        \item $\Tau$ acts on the stacks with no artificial drops;

        \item every normal tree on $(\M^+, \Sigma^+)$ of length $\leq \cont^+$ is also according to $\Tau$.
    \end{parts}

    \item Now, the tree $\T$ is according to $\Sigma$ and of length $< \cont^+$ (by Lemma \ref{2044}).
    This means that $\T$ is according to $\Tau$; so, $\Tau$ is as required.
\end{prooff}

\begin{corollary}\label{2022c}
    Let $M$ be super-suitable and let $\xi < \lh(\T_M) - 1$ be a limit ordinal.
    Then $\qq(\T_M \rest \xi)$ exists and is an initial segment of $\M^{\T_M}_\xi$.
\end{corollary}
\begin{proof}
    By Corollary \ref{2030}, $\delta(\T_M \rest \xi) < \Theta$.
    By Lemma \ref{2030'}, $\qq(\T_M \rest \xi, b)$ exists and is an initial segment of $\hod \para \Theta$.
    We will be done if we show that $\qq(\T_M \rest \xi, b)$ is countably iterable (this would mean that $\qq(\T_M \rest \xi) = \qq(\T_M \rest \xi, b)$).
    That $\qq(\T_M \rest \xi, b)$ is countably iterable will follow if we show that it is $(\cont^+ + 1)$-iterable.
    The latter is true by Proposition \ref{2011} because 
    $$\qq(\T_M \rest \xi, b) \unlhd \hod \para \Theta \unlhd \M_\infty.$$
\end{proof}


We characterized the branch $[0, \eta)_{\T_M}$ as the branch which has a countably iterable Q-structure.
If we put the tree $\T_M \rest \eta$ in a transitive model, take a hull, and collapse, then the hull of the tree should collapse to a countable tree on $M$ according to the short-tree strategy.\footnote{
The point is that the short-tree strategy choose branches according to (countably) iterable Q-structures.
}
To get that $[0, \eta)_{\T_M}$ belongs to $L(\R)$, the main idea is to notice that this process can be ``reversed'': we can reconstruct the tree based on its hulls.
The reconstruction is done by an ultrapower argument using the ultrafilter $\mu_\eta$.

We start by verifying that hulls of the tree collapse to a tree according to the short-tree strategy.

\begin{lemma}\label{2022a}
    Suppose
    \begin{assume}
        \item $M$ is super-suitable,
        
        \item $\xi < \lh(\T_M)$,
        
        \item $\theta$ is large enough regular,

        \item $X \prec (H_\theta, \in, \T_M \rest \xi)$ is countable,

        \item $\pi_X: \Bar{X} \to X$ is the anti-collapse.
    \end{assume}
    Then $\pi_X^{-1}(\T_M \rest \xi)$ is according to $\Sigma_M^\mathrm{st}$.
\end{lemma}
\begin{prooff}
    \item Assuming that $M$ is fixed, we reach the conclusion by induction on $\xi$.
    The only non-trivial case in this induction is when $\xi$ is a successor of a limit ordinal, so this is what we focus on below.

    \item\label{2047} Let $\T$ and $\alpha$ be the pre-images of $\T_M \rest \xi$ and $\xi - 1$ (resp.) by $\pi_X$.
    What needs to be verified is that $\qq(\T \rest \alpha)$ exists and is an initial segment of $\M^{\T}_\alpha$.
    More precisely, we want to find $Q$ such that
    \begin{parts}
        \item $\M(\T \rest \alpha) \unlhd Q$,
        
        \item $\delta(\T \rest \alpha)$ is a strong cutpoint in $Q$,
        
        \item $Q$ is the Q-structure of $Q$ at $\delta(\T \rest \alpha)$,

        \item $Q$ is sound above $\delta(\T \rest \alpha)$,

        \item \label{2059} $Q$ is iterable,

        \item $Q \unlhd \M^{\T}_\alpha$.
    \end{parts}

    \item By Lemma \ref{2022c}, $\qq(\T_M \rest (\xi - 1))$ exists and is an initial segment of $\M^{\T_M}_{\xi - 1}$.
    This means that $\qq(\T_M \rest (\xi - 1)) \in X$, so we may set $Q := \pi_X^{-1}(\qq(\T_M \rest \xi))$.

    \item Now, the facts that $\pi_X(\T) = \T_M \rest \xi$ and that $\pi_X$ is elementary, immediately imply all requirements from \ref{2047} except \ref{2059}.
    Requirement \ref{2059} follows from the facts that $\qq(\T_M \rest (\xi - 1))$ is countably iterable and that
    $$\pi_X \rest Q : Q \to \qq(\T_M \rest (\xi - 1))$$
    is elemenetary.
\end{prooff}

We now work towards articulating how to put together the branches through hulls of $\T_M \rest \eta$ to reconstruct the branch $[0, \eta)_{\T_M}$.
First of all, we will have to be careful how we code iteration trees.

\begin{remark}
    It will be relevant to know exactly what kind of information is coded by $\T$ when we say that $\T$ is a normal tree on a premouse $M$.
    We make the convention $\T$ in this case codes the 4-tuple
    $$(M,\ \lh(\T),\ (\lh(E^\T_\xi): \xi + 1 < \lh(\T)),\ ([0, \eta)_\T: \eta < \lh(\T) \mbox{ limit})).$$
    This information suffices to reconstruct all other information associated to $\T$ (the tree order, individual models, embeddings, etc.).
    \qed
\end{remark}

\begin{notation}
    We denote by $\langle \cdot, \cdot \rangle$ the Gödel's coding function $\ord \times \ord \to \ord.$
    \qed
\end{notation}

\begin{definition}
    Let $\T$ be a normal tree on a premouse $M$.
    We define the set of ordinals \intro{$\mathrm{code}(\T)$} as follows.
    \begin{parts}
        \item Let $<_M$ denote the constructibility order on $M$.
        For $\alpha, \beta < \otp(<_M)$, let us set $\alpha E \beta$ iff the $\alpha^\mathrm{th}$ element of $(M, <_M)$ belongs to the $\beta^\mathrm{th}$ element.
        Let $C_0 := \{\langle 0, \langle \alpha, \beta \rangle \rangle : \alpha E \beta\}$.

        \item Let $C_1 := \{\langle 1, \lh(\T) \rangle\}$.

        \item Let $C_2 := \{\langle 2, \lh(E^\T_\xi) \rangle: \xi + 1 < \lh(\T)\}$.

        \item Let $C_3 := \{ \langle 3, \langle \eta, \xi \rangle \rangle : \eta < \lh(\T)$ is limit and $\xi \in [0, \eta)_\T\}$.

        \item $\mathrm{code}(\T) := C_0 \cup C_1 \cup C_2 \cup C_3$.
        \qed
    \end{parts}
\end{definition}

\begin{remark}
    Let $\T$ be a normal tree on a premouse $M$.
    Then $\T \in J_1(\mathrm{code}(\T))$.
    On the other hand, letting $\alpha := \sup(\mathrm{code}(\T)) + 1$, we have that $\mathrm{code}(\T) \in J_\alpha (\T)$.
    \qed
\end{remark}

Now follows some relatively standard notation for hulls and their (anti-)collapses.

\begin{notation}
    Suppose that
    \begin{assume}
        \item $W\ni\eta$ is a rudimentarily closed transitive set,

        \item $\ww:=(W,\in,\dots)$ is a first-order structure in a countable language with a lightface definable wellordering,
        
        \item $\sigma\subseteq W$.
    \end{assume}
    Then we denote by
    \begin{parts}
        \item \intro{$\hull^\ww(\sigma)$} the Skolem hull of $\sigma$ inside $\ww$, understood both as a set and as a substructure,

        \item \intro{$\chull^\ww(\sigma)$} the transitive collapse of $\hull^\ww(\sigma)$,

        \item \intro{$\pi_\sigma^\ww$} the anti-collapse mapping
        $$\chull^\ww(\sigma) \longrightarrow \hull^\ww(\sigma).$$
        \qed
    \end{parts}
\end{notation}

Here is how to characterize the branch $[0, \eta)_{\T_M}$ using countable hulls of $\T_M \rest \eta$ and the short-tree strategy.

\begin{proposition}\label{2041}
    Suppose
    \begin{assume}
        \item $M$ is super-suitable,
        
        \item $\eta < \lh(\T_M) - 1$ is a limit ordinal,

        \item $\ww := (J_1(\mathrm{code}(\T_M \rest \eta)), \in, \mathrm{code}(\T_M \rest \eta))$,

        \item $\alpha < \eta$.
    \end{assume}
    Then $\alpha \in [0, \eta)_{\T_M}$ if and only if for club-many $\sigma \in [\eta]^\omega$, $(\pi^\ww_\sigma)^{-1}(\T_M \rest \eta) \in \dom(\Sigma_M^{\mathrm{st}})$ and $(\pi^\ww_\sigma)^{-1}(\alpha) \in \Sigma_M^{\mathrm{st}} \left( (\pi^\ww_\sigma)^{-1}(\T_M \rest \eta) \right)$.
\end{proposition}
\begin{prooff}
    \item Let $E$ consist of all $\sigma \in [\eta]^\omega$ such that $\alpha \in \sigma$.
    It is easily seen that $E$ contains a club.

    \item For $\sigma \in E$, let $U_\sigma$ and $\alpha_\sigma$ denote the pre-images of $\T_M \rest \eta$ and $\alpha$ (resp.) under $\pi^\ww_\sigma$.

    \item Let us assume that $\alpha \in [0, \eta)_{\T_M}$ and let us show that for club-many $\sigma \in E$, $U_\sigma \in \dom(\Sigma^\mathrm{st}_M)$ and $\alpha_\sigma \in \Sigma_M^\mathrm{st}(U_\sigma)$.

    \item Let $\theta$ be large enough regular.
    For $X \prec H_\theta$, let $\sigma_X := \eta \cap X$.
    We will achieve our goal if we show that for every countable
    $$X \prec (H_\theta, \in, \T_M \rest (\eta + 1), \alpha),$$
    $U_{\sigma_X} \in \dom(\Sigma^\mathrm{st}_M)$ and $\alpha_{\sigma_X} \in \Sigma_M^\mathrm{st}(U_{\sigma_X})$.

    \item Let us fix countable $X \prec (H_\theta, \in, \T_M \rest (\eta + 1), \alpha)$ and let us denote by $\pi_X: \Bar{X} \to X$ its anti-collapse.

    \item \claim $U_{\sigma_X} = \pi_X^{-1}(\T_M \rest \eta)$ and $\alpha_{\sigma_X} = \pi_X^{-1}(\alpha)$.
    \begin{prooff}
        \item Let $\beta := \sup(\mathrm{code}(\T_M \rest \eta)) + 1$.
        Observe that $\hull^\ww(\sigma) \cap \beta = X \cap \beta$.
        
        \item Let $\rho: \beta_X \to \beta$ be the anti-collapse of $X \cap \beta$.
        We conclude that
        $$\pi_X \rest \beta_X  = \pi^{\ww}_{\sigma_X} \rest \beta_X = \rho.$$

        \item Thus, $\alpha_{\sigma_X} = (\pi^\ww_{\sigma_X})^{-1}(\alpha) = \pi_X^{-1}(\alpha)$.

        \item We also conclude that for all $x \subseteq \beta_X$ with $x \in \chull^\ww(\sigma)$,
        $$\pi_X(x) = \pi^\ww_{\sigma_X}(x) = \rho [x].$$
        
        \item In particular,
        $$\pi_X (\mathrm{code}(U_{\sigma_X})) = \pi^\ww_{\sigma_X}(\mathrm{code}(U_{\sigma_X})) = \mathrm{code}(\T_M \rest \eta) = $$
        $$= \pi_X(\mathrm{code}(\pi_X^{-1}(\T_M \rest \eta))).$$

        \item In conclusion,
        $$\mathrm{code}(U_{\sigma_X}) = \mathrm{code}(\pi_X^{-1}(\T_M \rest \eta))$$
        and consequently, $U_{\sigma_X} = \pi_X^{-1}(\T_M \rest \eta)$.
    \end{prooff}

    \item The claim reduces our goal to showing that $\pi_X^{-1}(\T_M \rest \eta) \in \dom(\Sigma^\mathrm{st}_M)$ and $\pi_X^{-1}(\alpha) \in \Sigma_M^\mathrm{st} \left( \pi_X^{-1}(\T_M \rest \eta) \right)$.

    \item By Lemma \ref{2022a}, $\pi_X^{-1}(\T_M \rest (\eta + 1))$ is according to $\Sigma^\mathrm{st}_M$, so in particular, $\pi_X^{-1}(\T_M \rest \eta) \in \dom(\Sigma^\mathrm{st}_M)$ and 
    $$\Sigma_M^\mathrm{st} \left( \pi_X^{-1}(\T_M \rest \eta) \right) = \pi_X^{-1}([0, \eta)_{\T_M}) \ni \pi_X^{-1}(\alpha).$$

    \item With this, we have verified the implication (${\implies}$) of the lemma.
    In exactly the same way, we can show that if $\alpha \not\in [0, \eta)_{\T_M}$, then for club-many $\sigma \in [\eta]^\omega$, $U_\sigma \in \dom(\Sigma^\mathrm{st}_M)$ and $\alpha_\sigma \not\in \Sigma_M^\mathrm{st}(U_\sigma)$, which is more than enough for the implication (${\impliedby}$).
\end{prooff}

Finally, it is not hard to see that this characterization can be ``seen'' by $L(\R)$.

\begin{proposition}\label{2028}
    Let $M$ be super-suitable and let $\eta < \lh(\T_M)$ be a limit ordinal.
    Suppose $\T_M \rest \eta \in L(\R)$.
    Then $[0, \eta)_{\T_M} \in L(\R)$.
\end{proposition}
\begin{prooff}

    \item Let $\ww := (J_1(\mathrm{code}(\T_M \rest \eta)), \in, \mathrm{code}(\T_M \rest \eta))$.
    By the assumption of the lemma, $\ww \in L(\R)$.

    \item Recall that $([\eta]^\omega)^{L(\R)} = [\eta]^\omega$.\footnote{
    The point here is that $\eta^\omega \subseteq L(\R)$.
    To see this, fix $f \in L(\R)$ such that $f : \R \twoheadrightarrow \eta$.
    Since $\R^\omega \in L(\R)$, we conclude that
    $$F : \R^\omega \twoheadrightarrow \eta^\omega : \sigma \mapsto f \circ \sigma$$
    also belongs to $L(\R)$.
    Therefore, $\eta^\omega = \im(F) \in L(\R)$.
    }

    \item For $\sigma \in [\eta]^\omega$, let $U_\sigma$ be the pre-image of $\T_M \rest \eta$ under $\pi^\ww_\sigma$.
    Since $\T_M \rest \eta \in L(\R)$, we have that $U = (U_\sigma : \sigma \in [\eta]^\omega) \in L(\R)$.

    \item For $\alpha < \eta$, let $E_\alpha$ be the set of all $\sigma \in [\eta]^\omega$ such that $\alpha \in \hull^\ww(\sigma)$.
    It is easily seen that $E_\alpha$ is a club and that it belongs to $L(\R)$.
    Furthermore, $(E_\alpha: \alpha < \eta) \in L(\R)$.
    
    \item For $\alpha < \eta$ and $\sigma \in E_\alpha$, let $\alpha \proj \sigma$ denote the pre-image of $\alpha$ under $\pi^\ww_\sigma$.
    Note that $(\alpha \proj \sigma: \alpha < \eta, \sigma \in E_\alpha) \in L(\R)$.

    \item Notice that $\Sigma_M^{\mathrm{st}} \in L(\R)$.\footnote{
    The partial strategy $\Sigma_M^{\mathrm{st}}$ is essentially a set of reals, so it is a subset of $L(\R)$.
    To see that it is really an element of $L(\R)$, it suffices to verify that its definition is absolute between $V$ and $L(\R)$.
    This is true because, by Proposition \ref{Q-structures L(R)}, $L(\R)$ can identify the relevant Q-structures.
    }

    \item For $\alpha < \eta$, let
    $$X_\alpha := \{ \sigma \in E_\alpha: U_\sigma \in \dom(\Sigma^\mathrm{st}_M), \alpha \proj \sigma \in \Sigma^\mathrm{st}_M(U_\sigma)\}.$$
    Putting all of the above together, we conclude that for all $(X_\alpha: \alpha < \eta) \in L(\R)$.

    \item Finally, by Proposition \ref{2041}, we have that $\alpha \in [0, \eta)_{\T_M}$ if and only if $X_\alpha$ contains a club in $[\eta]^\omega$, which is true if and only if $X_\alpha \in \mu_\eta$.
    This shows that $[0, \eta)_{\T_M} \in L(\R)$, as required.
\end{prooff}

This proposition is the inductive step.
Before we put everything together into the final conclusion, recall that we also want to see that $\lh(\T_M) = \Theta + 1$.
The next lemma establishes the bound on how big models produced by iterations of a given size can be.

\begin{lemma}\label{2164}
    Assume only $\mathrm{ZF}$.
    Let $\mathcal{X}$ be a normal tree on a countable premouse with the last model $R$.
    Then $|R| \leq \max\{\lh(\mathcal{X}), \omega\}$.
\end{lemma}
\begin{proof}
    If $\mathrm{AC}$ holds, this easily follows by induction on $\lh(\mathcal{X})$.
    In general, notice that for some $X \subseteq \ord$, $\mathcal{X} \in L[X]$.
    Since $L[X] \models \mathrm{AC}$, we conclude that $L[X] \models |R| \leq \max\{\lh(\mathcal{X}), \omega\}$ and consequently, $V \models |R| \leq \max\{\lh(\mathcal{X}), \omega\}$.
\end{proof}

We now do one induction to show both that initial segments of the tree $\T_M$ belong to $L(\R)$ and that length of the tree is $\geq \Theta$.

\begin{proposition}\label{2035}
    Let $M$ be super-suitable.
    Then for all $\xi \leq \Theta$, $\lh(\T_M) \geq \xi$ and $\T_M \rest \xi \in L(\R)$.
\end{proposition}
\begin{prooff}
    \item We proceed by induction on $\xi$.
    The \textbf{case $\xi = 0$} is trivial.

    \item \textbf{Case where $\xi$ is limit.}
    The inductive hypothesis yields that $\lh(\T_M) \geq \xi$ and $\T_M \rest \xi \subseteq L(\R)$.


    \item Proposition \ref{2013} implies that $\T_M \rest \xi$ is a \textit{definable} subset of $L(\R)$.
    This suffices for the conclusion that $\T_M \rest \xi \in L(\R)$.

    \item \textbf{Case where $\xi$ is a successor of a limit ordinal.}
    The inductive hypothesis implies that $\lh(\T_M) \geq \xi - 1$ and $\T_M \rest (\xi - 1) \in L(\R)$.

    \item Since $\T_M$ has a last model, its length cannot be equal to the limit ordinal $\xi - 1$.
    Thus, $\lh(\T_M) \geq \xi$.

    \item By Proposition \ref{2028}, $[0, \xi - 1)_{\T_M} \in L(\R)$, so
    $$\T_M \rest \xi = (\T_M \rest (\xi - 1)) ^\frown [0, \xi - 1)_{\T_M} \in L(\R),$$
    as required.

    \item \textbf{Case where $\xi$ is a successor of a successor ordinal.}
    By the inductive hypothesis, $\lh(\T_M) \geq \xi - 1$ and $\T_M \rest (\xi - 1) \in L(\R)$.

    \item Since $\xi$ is a successor, we have that $\xi < \Theta$.
    Lemma \ref{2164} applied to $\T_M \rest (\xi - 1)$ in $L(\R)$ then implies that the last model of $\T_M \rest (\xi - 1)$ is not equal to $\hod \para \Theta$.
    This means that $\T_M \neq \T_M \rest (\xi - 1)$, so $\lh(\T_M) \geq \xi$.

    \item Now, $\T_M \rest \xi $ is obtained from $\T_M \rest (\xi - 1)$ by applying a single extender, so $\T_M \rest \xi \in L(\R)$, as required.

    \item The four case cover all possibilites, which concludes the proof.
\end{prooff}

We are now ready to state the main result of this section.

\begin{corollary}
    Let $M$ be super-suitable.
    Then 
    \begin{parts}
        \item $\lh(\T_M) = \Theta + 1$,

        \item $\T_M \rest \Theta \in L(\R)$,

        \item $\T_M \not \in L(\R)$.
    \end{parts}
\end{corollary}
\begin{proof}
    By Propositions \ref{2044} and \ref{2035}, $\lh(\T_M) \in \{\Theta, \Theta + 1\}$.
    Since $\T_M$ has a last model, we conclude that $\lh(\T_M) = \Theta + 1$.
    Proposition \ref{2035} also implies that $\T_M \rest \Theta \in L(\R)$.
    To see that the main branch of $\T_M$ does not belong to $L(\R)$, note that otherwise we would have the embedding
    $$\pi^{\T_M}:M\to\M_{\infty}$$
    which is continuous\footnote{\label{(15)}
    A $\Sigma_0$-embedding $i: M \to N$ between rudimentary closed, transitive sets/classes is said to be \textit{continuous} at an ordinal $\delta \in M$ just in case the mapping
    $$i \rest (M \cap \ord): M \cap \ord \to N \cap \ord$$
    is continuous at $\delta$, with respect to the order topology on both the domain and codomain.
    Considering that this restricted function is strictly increasing, the continuity at $\delta$ is equivalent to the assertion that $\sup(i[\delta]) = i(\delta)$.
    } at the Woodin and whose image would thus witness that
    $$\cof^{L(\R)}(\Theta)=\omega,$$
    whereas $\Theta$ is regular in $L(\R)$ (see \cite[Lemma 2.19]{koellner2010large}).
\end{proof}

\section{\texorpdfstring{``Localizable'' Definition of $\hod$ in $L(\R)$}{Localizable Definition of HOD in L(R)}}
\label{local definition of hod i}

In this section find a definition of $\hod$ in $L(\R)$ which uses only local information.
We call it ``localizable'' because it will take another section to check that it correctly relativizes to an initial segment of $L(\R)$.

Thus, we fix a $\hod$-cardinal $\eta$ and aim to describe $\hod\para\eta^{+\hod}$ from $\hod \para \eta$ using only the ``local'' information from the ambient universe $L(\R)$.
We will immediately restrict to $\eta$'s which are greater or equal to $\boldsymbol{\delta}^2_1$.
The reason for this is that our argument relies on short-tree iterability, for whose identification we need the sets of reals from $L_{\boldsymbol{\delta}^2_1}(\R)$.

\begin{declaration}\label{eta} 
    We fix a $\hod$-cardinal \intro{$\eta$} satisfying that $\boldsymbol{\delta}^2_1\leq\eta<\Theta$.\qed
\end{declaration}

\begin{notation}\label{1086}
    \hfill
    \begin{parts}
        \item \intro{$\hh_\eta$} denotes the premouse $\hod\para\eta$.



        \item \intro{$\qq_\eta$} denotes the Q-structure $\qq(\hod\para\Theta,\eta)$.
        \qed


    \end{parts}
\end{notation}

Starting from a super-suitable mouse $M$, we can reach $\hod \para \Theta$ by the iteration $\T_M$.
Thus, all models on a tail of this iteration agree with $\hod$ up to $\eta$.
It turns out that the tail actually begins with the index $\eta$ itself.

\begin{proposition}\label{2173}
    Suppose $M$ is super-suitable.
    Then $\eta$ is the least ordinal $i < \lh(\T_M)$ such that $\M^{\T_M}_i \para \eta = \hod \para \eta$.
\end{proposition}
\begin{prooff}
    \item Let $\gamma_M$ denote the least ordinal $i < \lh(\T_M) = \Theta + 1$ such that $\M^{\T_M}_i \para \eta = \hod \para \eta$.

    \item Let us assume towards a contradiction that $\gamma_M > \eta$.
    We will reach a contradiction by showing that $\lh(E_\eta^{\T_M}) \geq \eta$.\footnote{
    This is a contradiction because it means that the equality $\M^{\T_M}_i \para \eta = \hod \para \eta$ already occurs for $i = \eta < \gamma_M$.
    }


    \item Since $\T_M$ a normal tree, the lengths of the extenders used in it increase.
    This means that
    $$\lh(E^{\T_M}_\eta) \geq \sup \{\lh(E^{\T_M}_i): i < \eta\} \geq \eta.$$
    
    \item This is a contradiction, so it follows that $\gamma_M \leq \eta$.

    \item Let us now argue for the other inequality.
    Proposition \ref{2013} gives us that $\T_M \rest (\gamma_M + 1) \in L(M, \hod \para \tau)$, where $\tau < \Theta$ is an arbitrary cardinal of $\hod$ bigger than $o(\M^{\T_M}_{\gamma_M})$.

    \item Let $x$ be a real that codes $M$.
    It follows that $\T_M \rest (\gamma_M + 1) \in \hod_x$.

    \item By Lemma \ref{1269}, $\eta$ is a cardinal of $\hod_x$.
    
    \item Since $\M^{\T_M}_{\gamma_M} \para \eta = \hod \para \eta$, we have that $|\M^{\T_M}_{\gamma_M}|^{\hod_x} \geq \eta$.
    
    \item Lemma \ref{2164} applied in $\hod_x$ yields that $\gamma_M \geq \eta$, as required.
\end{prooff}

\begin{definition}\label{1315} 
    Suppose that $M$ is super-suitable.
    Then we define the following.
    \begin{parts}

        \item \intro{$\uu_M$} is the tree $\T_M\rest(\eta+1)$.
        
        \item \intro{$\pp_M$} is the premouse $\M^{\T_M}_{\eta}$.\qed
    \end{parts}
\end{definition}

The iteration $\uu_M$ is the part of the iteration $\T_M$ necessary to get $\hod\para \eta$; the rest of the iteration $\T_M$ happens above $\eta$.
We will see that if we are more careful in choosing the initial mouse $M$, then the last model $\pp_M$ of $\uu_M$ not only reconstructs $\hod \para \eta$, but even $\hod\para \eta^{+\hod}$.
This idea is at heart of the argument.

The more careful choice of $M$ that we mentioned consists in focusing on what we call \textit{$\eta$-exact} mice.


\begin{definition}\label{eta-exact}
    Suppose that $M$ is super-suitable.
    Then 
    \begin{parts}
        \item $M$ is \intro{$\eta$-exact} iff $\eta\in\ran(\pi_{M,\infty})$,

        \item if $M$ is $\eta$-exact, then \intro{$\eta_M$} denotes the preimage of $\eta$ by $\pi_{M,\infty}$.\qed
    \end{parts}
\end{definition}

The essence of the local definability is in the following proposition.

\begin{proposition}\label{gamma on the main branch} 
    Suppose that $M$ is super-suitable and $\eta$-exact.
    Then it holds that
    \begin{parts}
        \item\label{1319} $\eta\in b_M$,

        \item\label{1321} $\crit(\pi^{\T_M}_{\eta,\Theta})>\eta$,

        \item\label{1730} $\eta\in\ran(\pi^{\T_M}_{0,\eta})$,

        \item\label{1732} $\pp_M \para \eta^{+\pp_M} = \hod \para \eta^{+\hod}$.
    \end{parts}
\end{proposition}
\begin{prooff}
    \item Let us assume towards contradiction that \ref{1319} fails.
    Then there exist $\xi,\zeta\in b_M$ such that $\eta\in (\xi,\zeta)$ and $\xi$ is a $\T_M$-predecessor of $\zeta$.

    \item Note that $\zeta=\Bar{\zeta}+1$ for some $\Bar{\zeta}\geq\eta$ and 
    $$\M^{\T_M}_\zeta=\ult(\M^{\T_M}_\xi,E^{\T_M}_{\Bar{\zeta}}).$$

    \item We have that
    $$\crit(E^{\T_M}_{\Bar\zeta})<\lh(E^{\T_M}_\xi)\leq\eta\leq\lh(E^{\T_M}_{\Bar\zeta}),$$
    where the first inequality follows from the fact that $\T_M$ is normal, the second one follows by definition of $\eta$, and the third one follows from the fact $\M_{\Bar{\zeta}}^{\T_M}\para\eta=\M_{\eta}^{\T_M}\para\eta=\M_\infty\para\eta$.

    \item The previous two points then imply that $\eta\not\in\ran (\pi^{\T_M}_{\xi,\zeta})$, which contradicts the fact that $\eta\in\ran(\pi^{\T_M})$.

    \item Let us now verify part \ref{1321}.
    Let $E$ be the extender used at $\eta$ along $b_M$.
    We have that $\lh(E)\geq\eta$.

     \item The case $\crit(E)\leq\eta$ would imply $\eta\not\in\ran (\pi^{\T_M}_{\eta,\Theta})$, which is not possible since $\eta\in\ran(\pi^{\T_M})$.

    \item Thus, $\crit(\pi^{\T_M}_{\eta,\Theta})=\crit(E)>\eta$, which establishes part \ref{1321}.
    
    \item Parts \ref{1730} and \ref{1732} now easily follow.
\end{prooff}

The rest consists in approximating the above proposition, not only in $L(\R)$, but in its initial segment.
Recall first that $L(\R)$ does not have the relevant iteration strategies to identify super-suitable mice; thus, we replace that with suitable premice.
Once we do this, we are not able to make sense of $\uu_M$ as we did before.
Likely, $\uu_M$ has a different characterization: it is the $M$ side of the co-iteration between $M$ and $\hod \rest \Theta$, where the branches are chosen according to Q-structures.

\begin{definition} \label{2305'}
    Suppose that $M$ is a countable premouse.
    Then we define the following.
    \begin{parts}
        \item We say that ``\intro{$\uu_M$ converges}'' iff there exists a normal tree $\uu$ on $M$ of length $\eta + 1$ such that 
        \begin{parts}
            \item for all $i < \eta$, the least disagreement between $\M^\uu_i$ and $\hh_\eta$ is on the $\M^\uu_i$-side and $E^\uu_i$ is that disagreement,
    
            \item for all limit $i \leq \eta$, we have that the Q-structure $\qq(\uu\rest i,[0,i)_\uu)$ exists and is an initial segment of $\qq_\eta$,

            \item the branch $[0, \eta]_\uu$ does not drop and $\M^{\uu}_\eta \para \eta = \hh_\eta$.
        \end{parts}
        
        \item If it holds that ``$\uu_M$ converges'', then the tree \intro{$\uu_M$} on $M$ is the unique normal tree $\uu$ on $M$ of length $\eta + 1$ satisfying that
        \begin{parts}
            \item for all $i < \eta$, the least disagreement between $\M^\uu_i$ and $\hh_\eta$ is on the $\M^\uu_i$-side and $E^\uu_i$ is that disagreement,
    
            \item for all limit $i < \eta$, we have that the Q-structure $\qq(\uu\rest i,[0,i)_\uu)$ exists and is an initial segment of $\qq_\eta$.
        \end{parts}




        \item If $\uu_M$ converges, we denote by \dfn{$\pp_M$} the last model of $\uu_M$ and we denote by \dfn{$\sigma_M$}\mrg{$\pp_M$, $\sigma_M$} the mapping
        $$\pi^{\uu_M}:M\to\pp_M.$$



        \item The premouse $M$ is \intro{$\eta$-exact} iff $\uu_M$ converges, $\pp_M\rhd\hh_\eta$, and $\eta\in\ran(\sigma_M)$.

        \item\mrg{$\eta_M$} If $M$ is $\eta$-exact, we denote $\eta_M:=\sigma_M^{-1}(\eta)$.\qed
    \end{parts}
\end{definition}

\begin{remark}
    The tree $\uu_M$ is built by comparing $M$ to $\hh_\eta$, using the strategy given by Q-structures at limit stages.
    We want that $\hh_\eta$ does not move in this comparison and that $M$ comes on the winning side, which we ensure through the notion of convergence.
    If $M$ is super-suitable and $\eta$-exact in the sense of Definition \ref{eta-exact}, then it is $\eta$-exact in this new sense (see Proposition \ref{gamma on the main branch}) and the objects $\eta_M$, $\uu_M$, and $\pp_M$ correspond to those introduced in Definition \ref{eta-exact} and Definition \ref{1315}.
    In that case, the embedding $\sigma_M$ is exactly the embedding
    $$\pi^{\T_M}_{0,\eta}=\pi^{\uu_M}: M \to \pp_M.$$
    \qed
\end{remark}

We do not know if the rest of the argument can be realized by working with premice which are only suitable and $\eta$-exact.
Luckily, this does not matter for our purposes.
The premice we will look at are those which, in addition, have the property that all of their pseudo-iterates are also $\eta$-exact and, furthermore, have the same ``opinion'' on what $\hod \para \eta^+{\hod}$ should be.
(Reminder: $\eta$-exact super-suitable mice have the correct opinion on $\hod \para \eta^{+\hod}$, so in particular, all of them have the same opinion.)

\begin{definition}\label{eta-good}
    Suppose that $M$ is a countable premouse.
    Then $M$ is \intro{$\eta$-good} iff
    \begin{parts}
        \item $M$ is suitable,

        \item $M$ is short tree iterable,

        \item\label{1856} all normal pseudo-iterates $N$ of $M$ are $\eta$-exact and satisfy that
        $$\pp_N \para \eta^{+\pp_N} = \pp_M \para \eta^{+\pp_M}.$$
        \qed
    \end{parts}
\end{definition}

As expected, the following is true.

\begin{proposition}\label{super-suitable, very good} 
    Suppose that $M$ is super-suitable and $\eta$-exact.
    Then $M$ is $\eta$-good.
\end{proposition}
\begin{prooff}
    \item We already know that super-suitable premice are suitable and short tree iterable, so it remains to verify condition \ref{1856} of Definition \ref{eta-good}.

    \item Let $N$ be a normal pseudo-iterate of $M$.
    We want to show that $N$ is $\eta$-exact and that $\pp_N \para \eta^{+\pp_N} = \pp_M \para \eta^{+\pp_M}$.
    
    \item Let $\T$ be the tree satisfying $M\longrightarrow^\mathrm{st}_\T N$.
    By Lemma \ref{super-suitability and arrow}, we have that $\T$ is according to $\Sigma_M$ and $N$ is either the last model of $\T$ or $N=\M^{\T}_{\Sigma_M(\T)}$.
    In addition, the main branch of $\T$ or the branch $\Sigma_M(\T)$ does not drop.

    \item This implies that $N$ is super-suitable.
    Since $\eta\in\ran(\pi_{M,\infty})$ and $\pi_{M,\infty}=\pi_{N,\infty}\circ\pi_{M,N}$, we get that $\eta\in\ran(\pi_{N,\infty})$, which is another way of saying that $N$ is $\eta$-exact.

    \item By Proposition \ref{gamma on the main branch}, we have that
    $$\pp_N \para \eta^{+\pp_N} = \hod \para \eta^{+\hod} =\pp_M \para \eta^{+\pp_M},$$
    which concludes the verification.
\end{prooff}

Here is a characterization of $\hod \para \eta^{+\hod}$ inside $L(\R)$ which will turn out to be ``local enough''.

\begin{proposition}\label{characterization}
    Let $M$ be an $\eta$-good premouse.
    Then $\pp_M \para \eta^{+\pp_M} = \hod \para \eta^{+\hod}$.
\end{proposition}
\begin{prooff}


    \item There exists a super-suitable $N$ such that $\eta\in\ran(\pi_{N,\infty})$.\footnote{\label{(19)}
    To see this, recall that $\eta < \Theta$, as assumed in Declaration \ref{eta}.
    By Theorem \ref{1639}, $\eta \in \M_\infty$, so by Definition \ref{1630'} of $\M_\infty$, there are $N \in \ii(\M)$ and $\Bar{\eta} \in N$ such that $\pi_{N, \infty}(\Bar{\eta}) = \eta$.
    Note that $N$ is as required (recall Definition \ref{724} of being super-suitable).
    }

    \item By Proposition \ref{comparison}, there exists a suitable $P$ which is a normal pseudo-iterate of both $M$ and $N$.

    \item\label{1529} Since $M$ is $\eta$-good, we have that
    $$\pp_M \para \eta^{+\pp_M} = \pp_P \para \eta^{+\pp_P}.$$

    \item Since $N$ is super-suitable and $\eta\in\ran(\pi_{N,\infty})$, we have that $P$ is super-suitable as well and $\eta\in\ran(\pi_{P,\infty})$.

    \item\label{1534} By Proposition \ref{gamma on the main branch}, we have that 
    $$\pp_P \para \eta^{+\pp_P} = \hod \para \eta^{+\hod}.$$

    \item By \ref{1529} and \ref{1534}, we have that
    $$\pp_M \para \eta^{+\pp_M} = \hod \para \eta^{+\hod}.$$
\end{prooff}

\section{\texorpdfstring{Local Definition of $\hod$ in $L(\R)$}{Local Definition of HOD in L(R)}}
\label{local definition of hod}

In the previous section, we reduced the problem to examining to what extent the definition of ``$\eta$-good'' is local.
By looking at it more carefully, we see that, in order to make sense of it, we essentially need only two things: we need short-tree strategies and we need $\qq\eta$.
It is now time to make explicit in which structure we want to define $\hod\para\eta^{+\hod}$.

\begin{definition}\label{1359}
    We use the following notation.
    \begin{assume}

        \item \intro{$\chi_\eta$} denotes the least ordinal $\chi>\eta$ satisfying that 
        $$L_\chi (\hh_\eta^\omega)[\mu_\eta]\models\mathrm{ZF}^-+\mbox{``}\ps([\eta]^\omega)\mbox{ exists''}$$
        and $\mu_\eta \cap L_\chi (\hh_\eta^\omega)[\mu_\eta]$ belongs to $L_\chi (\hh_\eta^\omega)[\mu_\eta]$ and is a supercompactness measure there.

        \item We write \intro{$L_{\chi_\eta}(\hh_\eta^\omega)[\mu_\eta]$} to mean the structure 
        $$(L_{\chi_\eta}(\hh_\eta^\omega)[\mu_\eta],\in,\mu_\eta,\hh_\eta),$$
        where (the restriction of) $\in$ is a binary predicate, (the restriction of) $\mu_\eta$ is a unary predicate, and $\hh_\eta$ is a constant.\qed
    \end{assume}
\end{definition}

\begin{lemma}\label{2151}
    The ordinal $\chi_\eta$ is well defined and strictly less than $\Theta$.
\end{lemma}
\begin{prooff}
    \item We work in $L(\R)$ and we denote by $\theta$ the cardinal $\Theta^{+10}$.

    \item\label{1666} Let us first verify that $L_\theta (\R)$ satisfies $\mathrm{ZF}^-$.
    What is nontrivial is to see that this model satisfies Collection, so we concentrate on that.

    \item Let $x\in L_\theta (\R)$ and $R\subseteq L_\theta (\R)^2$ be arbitrary.
    We want to find $y\in L_\theta (\R)$ such that for all $u\in x\cap\dom(R)$, $R[u]\cap y\not=\emptyset$.


    \item For all $u\in x$, let $g(u)$ be the least $\eta<\theta$ such that $R[u]\cap L_\eta(\R)\not=\emptyset$.
    We have defined a function $g:x\to \theta$.

    \item\claim $g$ is not cofinal.
    \begin{prooff}
        \item Let us assume otherwise.
        
        \item There exist $\xi<\theta$ and a surjection $f:\xi\times\R\twoheadrightarrow x$.

        \item Let $h:\omega\twoheadrightarrow\R$ code a generic for $\col(\omega,\R)$ and let us work in $L(\R)[h]=L[h]$.
        We have that $\theta$ remains a successor.
        
        \item Since Choice holds, we get that $\theta$ regular.

        \item The mapping $g\circ f\circ (\id_\xi\times h):\xi\times\omega\to \theta$ is cofinal.
        This contradicts the regularity of $\theta$.
    \end{prooff}

    \item Let $\zeta:=\sup(\ran(g))<\theta$ and let $y:=L_\zeta(\R)$.
    It is immediate that $y$ is as required in \ref{1666}.

    \item We have established that $L_\theta(\R)\models\mathrm{ZF}^-$.
    By Theorem \ref{1896},\footnote{
    Notice that $\hh_\eta^\omega \subseteq L(\R)$.
    To see this, recall that $\hh_\eta$ is a premouse, so the conclusion follows from the fact that $\eta^\omega \subseteq L(\R)$ (proved in Lemma \ref{1889}).
    } we have
    $$L_\theta (\hh_\eta^\omega)[\mu_\eta] = L_\theta (\R) \models \mathrm{ZF}^-+\mbox{``}\ps([\eta]^\omega)\mbox{ exists.''}$$
    
    \item Since there exists a surjection $\theta\times\R\twoheadrightarrow L_\theta(\R)$ and there is a surjection $\R \twoheadrightarrow \hh_\eta^\omega$,\footnote{
    Since $\hh_\eta$ is a premouse, there is a simply definable surjection of $\eta$ onto $\hh_\eta$.
    Now, the situation is as follows: there is a bijection from $\R$ onto $\R^\omega$; there is a surjection of $\R^\omega$ onto $\eta^\omega$; there is a surjection of $\eta^\omega$ onto $\hh_\eta^\omega$.
    } there is an elementary submodel
    $$X\prec L_\theta (\hh_\eta^\omega)[\mu_\eta]$$ 
    such that:
    \begin{itemize}
        \item $\hh_\eta^\omega \subseteq X$;

        \item $\mu_\eta \cap X \in X$;

        \item there is a surjection $\R \twoheadrightarrow X$.
    \end{itemize}

    \item The transitive collapse of $X$ is of the form
    $$L_\chi (\hh_\eta^\omega)[\mu_\eta]$$
    with $\chi<\Theta$.
    This suffices for the conclusion.
\end{prooff}

The first step is to check that $L_{\chi_\eta}(\hh_\eta^\omega)[\mu_\eta]$ has short-tree strategies.

\begin{lemma}\label{1516}
    The function $\lp$ is lightface definable over $L_{\chi_\eta} (\hh_\eta^\omega)[\mu_\eta]$.
\end{lemma}
\begin{proof}
    Definition \ref{Sigma-1 reflection} implies that $\Sigma_1$-formulas with parameters in $\R\cup\{\R\}$ are absolute between $L_{\chi_\eta} (\hh_\eta^\omega)[\mu_\eta]$ and $L(\R)$, from which it easily follows that the definition of $\lp$ is absolute between $L_{\chi_\eta} (\hh_\eta^\omega)[\mu_\eta]$ and $L(\R)$.
\end{proof}

\begin{corollary}\label{3229}
    For every super-suitable $M$, we have $\Sigma^\mathrm{st}_M \in L_{\chi_\eta} (\hh_\eta^\omega)[\mu_\eta]$.
    \qed
\end{corollary}

The next step is to see that $L_{\chi_\eta}(\hh_\eta^\omega)[\mu_\eta]$ contains $\qq_\eta$ as an element.

\begin{proposition}
    $\qq_\eta \in L_{\chi_\eta} (\hh_\eta^\omega)[\mu_\eta]$
\end{proposition}
\begin{proof}
    Let us fix an $\eta$-exact, super-suitable premouse $M$.
    By induction on $\xi \leq (\eta + 1)$, we show that $\T_M \rest \xi \in L_{\chi_\eta} (\hh_\eta^\omega)[\mu_\eta]$.
    The only non-trivial step is when $\xi = \gamma + 1$, where $\gamma$ is a limit ordinal.
    This is by the proof of Proposition \ref{2028}, which we reproduce below.
    
    \begin{pfenum}
        \item Let $\thickmuskip=5mu plus 5mu minus 1mu \ww := (J_1(\mathrm{code}(\T_M \rest \gamma)), \in, \mathrm{code}(\T_M \rest \gamma))$.
        We have $\ww \in L_{\chi_\eta} (\hh_\eta^\omega)[\mu_\eta]$.
    
    
        \item For $\sigma \in [\eta]^\omega$, let $U_\gamma$ be the pre-image of $\T_M \rest \gamma$ under $\pi^\ww_\sigma$.
        Since $\T_M \rest \gamma \in L(\R)$, we have that $U = (U_\sigma : \sigma \in [\gamma]^\omega) \in L_{\chi_\eta} (\hh_\eta^\omega)[\mu_\eta]$.
    
        \item For $\alpha < \gamma$, let $E_\alpha$ be the set of all $\sigma \in [\gamma]^\omega$ such that $\alpha \in \hull^\ww(\sigma)$.
        It is easily seen that $E_\alpha$ is a club and that it belongs to $L_{\chi_\eta} (\hh_\eta^\omega)[\mu_\eta]$.
        Furthermore, $(E_\alpha: \alpha < \gamma) \in L_{\chi_\eta} (\hh_\eta^\omega)[\mu_\eta]$.
        
        \item For $\alpha < \gamma$ and $\sigma \in E_\alpha$, let $\alpha \proj \sigma$ denote the pre-image of $\alpha$ under $\pi^\ww_\sigma$.
        Note that $(\alpha \proj \sigma: \alpha < \gamma, \sigma \in E_\alpha) \in L_{\chi_\eta} (\hh_\eta^\omega)[\mu_\eta]$.
    
        \item By Corollary \ref{3229}, we have that $\Sigma_M^{\mathrm{st}} \in L_{\chi_\eta} (\hh_\eta^\omega)[\mu_\eta]$.
    
        \item For $\alpha < \gamma$, let
        $$X_\alpha := \{ \sigma \in E_\alpha: U_\sigma \in \dom(\Sigma^\mathrm{st}_M), \alpha \proj \sigma \in \Sigma^\mathrm{st}_M(U_\sigma)\}.$$
        Putting all of the above together, we conclude that for all $(X_\alpha: \alpha < \gamma) \in L_{\chi_\eta} (\hh_\eta^\omega)[\mu_\eta]$.
    
        \item Finally, by Proposition \ref{2041}, we have that $\alpha \in [0, \gamma)_{\T_M}$ if and only if $X_\alpha$ contains a club in $[\gamma]^\omega$, which is true if and only if $X_\alpha \in \mu_\gamma$.
        This shows that $[0, \gamma)_{\T_M} \in L_{\chi_\eta} (\hh_\eta^\omega)[\mu_\eta]$, as required.
    \end{pfenum}

    This finishes the induction.
    Thus, $\T_M \rest (\eta + 1) \in L_{\chi_\eta} (\hh_\eta^\omega)[\mu_\eta]$.
    Since $\qq_\eta$ is easily computed from $\M^{\T_M}_\eta$, we conclude that $\qq_\eta \in L_{\chi_\eta} (\hh_\eta^\omega)[\mu_\eta]$.
\end{proof}

Based on this, we are very closed to the conclusion, but we need to use $\qq_\eta$ as an additional parameter.

\begin{corollary}\label{3142}
    The following sets are contained in $L_{\chi_\eta}(\hh_\eta^\omega)[\mu_\eta]$ and definable over it, with the additional parameter $\qq_\eta$:
    \begin{parts}
        \item $\{M : M\mbox{ is suitable}\}$

        \item $\{M : M\mbox{ is short tree iterable}\}$,

        \item $\{M : M\mbox{ is }\eta\mbox{-good}\}$.
        \qed
    \end{parts}
\end{corollary}

\begin{corollary} 
    The premouse $\hod\para\eta^{+\hod}$ belongs to $L_{\chi_\eta}(\hh_\eta^\omega)[\mu_\eta]$ and definable over it from the parameter $\qq_\eta$.
\end{corollary}
\begin{proof}
    By Proposition \ref{characterization}, the premouse $\hod\para\eta^{+\hod}$ is the common value of all $\pp_M \para \eta^{+\pp_M}$, for $\eta$-good premice $M$.
    By Corollary \ref{3142}, this characterization relativizes to $L_{\chi_\eta} (\hh_\eta^\omega)[\mu_\eta]$.
\end{proof}

But actually, it is easy to see that $\qq_\eta$ is unnecessary.

\begin{proposition}
    The premouse $\qq_\eta$ is definable over the structure $L_{\chi_\eta} (\hh_\eta^\omega)[\mu_\eta]$ (with no additional parameters).
\end{proposition}
\begin{proof}
    Notice that $\qq_\eta$ is the Q-structure above $\hh_\eta$ (in the sense of Definition \ref{Q-structure above}).
    Since $L_{\chi_\eta} (\hh_\eta^\omega)[\mu_\eta]$ has the strategies for the hulls of this structure (see Proposition \ref{Q-structures L(R)}), this characterization relativizes correctly to $L_{\chi_\eta} (\hh_\eta^\omega)[\mu_\eta]$.
\end{proof}

Therefore, we get the main theorem.

\begin{corollary}\label{2715}
\label{main}
    The premouse $\hod\para\eta^{+\hod}$ is definable over the structure $L_{\chi_\eta}(\hh_\eta^\omega)[\mu_\eta]$ with no additional parameters.
    \qed 
\end{corollary}

Here is one application of the theorem.

\begin{corollary}
    $\eta^{+\hod}<\chi_\eta$.\qed
\end{corollary}

\newpage

\printindex
\end{document}